# REVERSIBILITY OF CHORDAL SLE

By Dapeng Zhan

*University of California, Berkeley*

We prove that the chordal $SLE_\kappa$ trace is reversible for $\kappa \in (0,4]$.

**1. Introduction.** Stochastic Loewner evolutions (SLEs) are introduced by Oded Schramm [11] to describe the scaling limits of some lattice models, whose scaling limits satisfy conformal invariance and Markov property. The basic properties of SLE are studied in [9]. There are several different versions of SLE. A chordal SLE is defined in a simply connected domain, which is about some random curve in the domain that grows from one boundary point to another.

So far it has been proved that the chordal $SLE_6$ is the scaling limit of the explorer line of the site percolation on the triangular lattice with half open and half closed boundary conditions ([13] and [2]); the chordal $SLE_8$ is the scaling limit of UST Peano curve with half free and half wired boundary conditions [6]; the chordal $SLE_4$ is the scaling limit of the contour line of the two-dimensional discrete Gaussian free field with suitable boundary values [12]; and the chordal $SLE_2$ is the scaling limit of LERW started near one boundary point, conditioned to leave the domain near the other boundary point [17]. In [5], the $SLE_{8/3}$ is proved to satisfy the restriction property. From these results, we know that the chordal $SLE_\kappa$ trace is reversible for $\kappa = 6, 8, 4, 2, 8/3$.

In [9], it is conjectured that the chordal $SLE_\kappa$ trace is reversible for all $\kappa \in [0, 8]$. Scott Sheffield proposed that the reversibility can be derived from the relationship with the Gaussian free field [10]. In this paper we will prove this conjecture for $\kappa \in (0, 4]$ using only techniques of probability theory and stochastic processes. The main idea of this paper is as follows.

Suppose $(\beta(t))$ is a chordal $SLE_\kappa$ trace in a simply connected domain $D$ from a prime end $a$ to another prime end $b$. From the Markov property of SLE, for a fixed time $t_0$, conditioned on the curve $\beta([0, t_0])$, the rest of the curve $(\beta(t)\colon t\geq t_0)$ has the same distribution as a chordal $\mathrm{SLE}_\kappa$ trace in $D_{t_0}:=D\setminus\beta([0,t_0])$ from $\beta(t_0)$ to $b$. Assume that the chordal $\mathrm{SLE}_\kappa$ trace is reversible. Then the reversal of $(\beta(t):t\geq t_0)$ has the same distribution as the chordal $\mathrm{SLE}_\kappa$ trace in $D_{t_0}$ from $b$ to $\beta(t_0)$. On the other hand, since $(\beta(t)\colon t\geq t_0)$ is a part of the $\mathrm{SLE}_\kappa$ trace in $D$ from $a$ to $b$, so from the reversibility, the reversal of $(\beta(t)\colon t\geq t_0)$ should be a part of $\mathrm{SLE}_\kappa$ trace in $D$ from $b$ to $a$. Suppose $\gamma$ is an $\mathrm{SLE}_\kappa$ trace in $D_{t_0}$ from $b$ to $\beta(t_0)$. From the above discussion, if we integrate $\gamma$ against all possible curves $\beta([0,t_0])$, we should get a part of the $\mathrm{SLE}_\kappa$ trace in $D$ from $b$ to $a$, assuming that the chordal $\mathrm{SLE}_\kappa$ trace is reversible.

To prove the reversibility, we want to find a coupling of two $\mathrm{SLE}_\kappa$ traces, one is from $a$ to $b$, the other is from $b$ to $a$, such that the two curves visit the same set of points. If such coupling exists, we choose a pair of disjoint hulls, each of which contains some neighborhood of $a$ or $b$ in $\mathbb{H}$, and stop the two traces when they leave one of the two hulls, respectively. Before these stopping times, the two traces are disjoint from each other. The joint distribution of the two traces up to these stopping time should agree with that of $\beta$ and $\gamma$ discussed in the last paragraph up to the same stopping times. The Girsanov Theorem suggests that this distribution is absolutely continuous w.r.t. that of two independent chordal $\mathrm{SLE}_\kappa$ traces (one from $a$ to $b$, the other from $b$ to $a$) stopped on leaving the above two hulls. And the Radon–Nikodym derivative is described by a two-dimensional local martingale, which has the property that when one variable is fixed, it is a local martingale in the other variable. This is the $M(\cdot,\cdot)$ in Theorem 4.1. It is closely related with Julien Dubédat's work about commutation relations for SLE [3].

Using the $M(\cdot,\cdot)$, we may construct a portion of the coupling up to certain stopping times. To construct the global coupling, the difficulty arises when the two hulls collide, and the absolute continuity blows up after that time. In fact, we can not expect that the global coupling we are looking for is absolutely continuous w.r.t. two independent SLE. Instead, the coupling measure will be the weak limit of a sequence of absolutely continuous coupling measures. Each measure in the sequence is generated from some two-dimensional bounded martingale, which is the $M_*(\cdot,\cdot)$ in Theorem 6.1. The important property of $M_*$ is that, on the one hand, it carries the information of $M$ as much as we want; on the other hand, it is uniformly bounded, and remains to be a martingale even after the two hulls collide. So $M_*$ can be used as the Radon–Nikodym derivative to define a global coupling measure.

It is known that, for $\kappa>8$, the chordal $\mathrm{SLE}_\kappa$ trace is not reversible [15]. So far the reversibility for $\kappa\in(4,8)$ is still unknown. Although the results about martingales in this paper hold for all $\kappa>0$, the argument in the last step of the proof essentially uses the property that, for $\kappa\in(0,4]$, the chordal $\text{SLE}_\kappa$ trace does not touch the boundary other than the initial and end points.

The technique developed in this paper may have other usage. For example, it is used in [15] to prove the Duplantier's duality conjecture about SLE. It may also be used to study the reversal of the trace of other variations of SLE, for example, $\text{SLE}(\kappa,\rho)$ [5], continuous LERW [14] and annulus SLE [15, 16].

This paper is organized in the following way. In Section 2 we give the definition of the chordal SLE and some other basic notation, and then present the main theorem of this paper. In Section 3 we study the relations of two SLE that grow in the same domain. In Section 4 we present the two-dimensional local martingale $M$, and check its property by direct calculation of stochastic analysis. In Section 5 we give some stopping times up to which $M$ is bounded. And at the end of Section 5 we give a detailed explanation of the meaning of $M$. In Section 6 we use the local martingale to construct some two-dimensional bounded martingale $M_*$. In Section 7 we use $M_*$ to construct a sequence of coupling measures. The limit of these measures in some suitable sense is also a coupling measure. We finally prove that, under the limit measure, the two $\text{SLE}_\kappa$ traces coincide with each other.

**2. Chordal Loewner equation and chordal SLE.** Let $\mathbb{H} = \{z \in \mathbb{C} : \operatorname{Im} z > 0\}$ denote the upper half complex plane. If $H$ is a bounded closed subset of $\mathbb{H}$ such that $\mathbb{H} \setminus H$ is simply connected, then we call $H$ a hull in $\mathbb{H}$ w.r.t. $\infty$. For such $H$, there is a unique $\varphi_H$ that maps $\mathbb{H} \setminus H$ conformally onto $\mathbb{H}$ such that $\varphi_H(z) = z + \frac{c}{z} + O(1/z^2)$ as $z \to \infty$ for some $c \geq 0$. Such $c$ is called the half-plane capacity of $H$, and is denoted by $\operatorname{hcap}(H)$.

PROPOSITION 2.1. *Suppose $\Omega$ is an open neighborhood of $x_0 \in \mathbb{R}$ in $\mathbb{H}$. Suppose $W$ maps $\Omega$ conformally into $\mathbb{H}$ such that, for some $r > 0$, if $z \to (x_0 - r, x_0 + r)$ in $\Omega$, then $W(z) \to \mathbb{R}$. So $W$ extends conformally across $(x_0 - r, x_0 + r)$ by the Schwarz reflection principle. Then for any $\varepsilon > 0$, there is some $\delta > 0$ such that if a hull $H$ in $\mathbb{H}$ w.r.t. $\infty$ is contained in $\{z \in \mathbb{H} : |z - x_0| < \delta\}$, then $W(H)$ is also a hull in $\mathbb{H}$ w.r.t. $\infty$, and*

$$|\operatorname{hcap}(W(H)) - W'(x_0)^2 \operatorname{hcap}(H)| \leq \varepsilon |\operatorname{hcap}(H)|.$$

PROOF. This is Lemma 2.8 in [4]. □

For a real interval $I$, let $C(I)$ denote the real-valued continuous function on $I$. Suppose $\xi \in C([0,T))$ for some $T \in (0,+\infty]$. The chordal Loewner equation driven by $\xi$ is as follows:

$$(2.1) \qquad \partial_t \varphi(t,z) = \frac{2}{\varphi(t,z) - \xi(t)}, \qquad \varphi(0,z) = z.$$

For $0 \le t < T$, let $K(t)$ be the set of $z \in \mathbb{H}$ such that the solution $\varphi(s, z)$ blows up before or at time $t$. We call $K(t)$ and $\varphi(t, \cdot)$, $0 \le t < T$, chordal Loewner hulls and maps, respectively, driven by $\xi$. Then for each $t \in [0, T)$, $\varphi(t, \cdot)$ maps $\mathbb{H} \setminus K(t)$ conformally onto $\mathbb{H}$. Suppose for every $t \in [0, T)$,

$$\beta(t) := \lim_{z \in \mathbb{H}, z \to \xi(t)} \varphi(t, \cdot)^{-1}(z) \in \mathbb{H} \cup \mathbb{R}$$

exists, and $\beta(t)$, $0 \le t < T$, is a continuous curve. Then for every $t \in [0, T)$, $K(t)$ is the complement of the unbounded component of $\mathbb{H} \setminus \beta((0, t])$. We call $\beta$ the chordal Loewner trace driven by $\xi$. In general, such trace may not exist.

We say $(K(t), 0 \le t < T)$ is a Loewner chain in $\mathbb{H}$ w.r.t. $\infty$, if each $K(t)$ is a hull in $\mathbb{H}$ w.r.t. $\infty$; $K(0) = \varnothing$; $K(s) \subsetneq K(t)$ if $s < t$; and for each fixed $a \in (0, T)$, the extremal length [1] of the curve in $\mathbb{H} \setminus K(t+\varepsilon)$ that disconnects $K(t+\varepsilon) \setminus K(t)$ from $\infty$ tends to $0$ as $\varepsilon \to 0^+$, uniformly in $t \in [0, a]$. If $u(t)$, $0 \le t < T$, is a continuous (strictly) increasing function, and satisfies $u(0) = 0$, then $(K(u^{-1}(t)), 0 \le t < u(T))$ is also a Loewner chain in $\mathbb{H}$ w.r.t. $\infty$, where $u(T) := \sup u([0, T))$. It is called the time-change of $(K(t))$ through $u$. Here is a simple example of the Loewner chain. Suppose $\beta(t)$, $0 \le t < T$, is a simple curve with $\beta(0) \in \mathbb{R}$ and $\beta(t) \in \mathbb{H}$ for $t \in (0, T)$. Let $K(t) = \beta((0, t])$ for $0 \le t < T$. Then $(K(t), 0 \le t < T)$ is a Loewner chain in $\mathbb{H}$ w.r.t. $\infty$. It is called the Loewner chain generated by $\beta$.

If $H_1 \subset H_2$ are two hulls in $\mathbb{H}$ w.r.t. $\infty$, let $H_2/H_1 := \varphi_{H_1}(H_2 \setminus H_1)$. Then $H_2/H_1$ is also a hull in $\mathbb{H}$ w.r.t. $\infty$, $\varphi_{H_2/H_1} = \varphi_{H_2} \circ \varphi_{H_1}^{-1}$, and $\mathrm{hcap}(H_2/H_1) = \mathrm{hcap}(H_2) - \mathrm{hcap}(H_1)$. If $H_1 \subset H_2 \subset H_3$ are three hulls in $\mathbb{H}$ w.r.t. $\infty$, then $H_2/H_1 \subset H_3/H_1$ and $(H_3/H_1)/(H_2/H_1) = H_3/H_2$.

PROPOSITION 2.2. (a) *Suppose $K(t)$ and $\varphi(t, \cdot)$, $0 \le t < T$, are chordal Loewner hulls and maps, respectively, driven by $\xi \in C([0, T))$. Then $(K(t), 0 \le t < T)$ is a Loewner chain in $\mathbb{H}$ w.r.t. $\infty$, $\varphi_{K(t)} = \varphi(t, \cdot)$, and $\mathrm{hcap}(K(t)) = 2t$ for any $0 \le t < T$. Moreover, for every $t \in [0, T)$,*

(2.2) $$\{\xi(t)\} = \bigcap_{\varepsilon \in (0, T-t)} \overline{K(t+\varepsilon)/K(t)}.$$

(b) *Let $(L(s), 0 \le s < S)$ be a Loewner chain in $\mathbb{H}$ w.r.t. $\infty$. Let $v(s) = \mathrm{hcap}(L(s))/2$, $0 \le s < S$. Then $v$ is a continuous increasing function with $u(0) = 0$. Let $T = v(S)$ and $K(t) = L(v^{-1}(t))$, $0 \le t < T$. Then $K(t)$, $0 \le t < T$, are chordal Loewner hulls driven by some $\xi \in C([0, T))$.*

PROOF. This is almost the same as Theorem 2.6 in [4]. □

Let $B(t)$ be a (standard linear) Brownian motion, $\kappa \in (0, \infty)$, and $\xi(t) = \sqrt{\kappa} B(t)$, $0 \le t < \infty$. Let $K(t)$ and $\varphi(t, \cdot)$, $0 \le t < \infty$, be the chordal Loewner hulls and maps, respectively, driven by $\xi$. Then we call $(K(t))$ the standard chordal $\text{SLE}_\kappa$. From [9], the chordal Loewner trace $\beta(t)$, $0 \leq t < \infty$, driven by $\xi$ exists almost surely. Such $\beta$ is called the standard chordal $\text{SLE}_\kappa$ trace. We have $\beta(0) = 0$ and $\lim_{t\to\infty} \beta(t) = \infty$. If $\kappa \in (0,4]$, then $\beta$ is a simple curve, $\beta(t) \in \mathbb{H}$ for $t > 0$, and $K(t) = \beta((0,t])$ for $t \geq 0$; if $\kappa \in (4,\infty)$, then $\beta$ is not a simple curve. If $\kappa \in [8,\infty)$, then $\beta$ visits every $z \in \overline{\mathbb{H}}$; if $\kappa \in (0,8)$, then the Lebesgue measure of the image of $\beta$ in $\mathbb{C}$ is 0.

Suppose $D \subsetneq \mathbb{C}$ is a simply connected domain, and $a \neq b$ are two prime ends [1] of $D$. Then there is $W$ that maps $(\mathbb{H}; 0, \infty)$ conformally onto $(D; a, b)$. We call the image of the standard chordal $\text{SLE}_\kappa$ under $W$ the chordal $\text{SLE}_\kappa$ in $D$ from $a$ to $b$, which is denoted by $\text{SLE}_\kappa(D; a \to b)$. Such $W$ is not unique, but the $\text{SLE}_\kappa(D; a \to b)$ defined through different $W$ have the same distribution up to a linear time-change because the standard chordal $\text{SLE}_\kappa$ satisfies the scaling property. The main theorem in this paper is as follows.

THEOREM 2.1. *Suppose $\kappa \in (0,4]$, $\beta_1(t)$, $0 \leq t < \infty$, is an $\text{SLE}_\kappa(D; a \to b)$ trace, and $\beta_2(t)$, $0 \leq t < \infty$, is an $\text{SLE}_\kappa(D; b \to a)$ trace. Then the set $\{\beta_1(t): 0 < t < \infty\}$ has the same distribution as $\{\beta_2(t): 0 < t < \infty\}$.*

**3. Ensemble of two chordal Loewner chains.** In this section we study the relations of two chordal Loewner chains that grow together. Some computations were done in [3, 4, 5] and other papers. We will give self-contained arguments for all results in this section. Suppose $K_j(t)$ and $\varphi_j(t,\cdot)$, $0 \leq t < S_j$, are chordal Loewner hulls and maps, respectively driven by $\xi_j \in C([0, S_j))$, $j = 1, 2$. Assume that for any $t_1 \in [0, S_1)$ and $t_2 \in [0, S_2)$, $\overline{K_1(t_1)} \cap \overline{K_2(t_2)} = \varnothing$, then $K_1(t_1) \cup K_2(t_2)$ is a hull in $\mathbb{H}$ w.r.t. $\infty$. Fix $j \neq k \in \{1, 2\}$ and $t_0 \in [0, S_k)$. For $0 \leq t < S_j$, let

(3.1) $\qquad K_{j,t_0}(t) = (K_j(t) \cup K_k(t_0))/K_k(t_0) = \varphi_k(t_0, K_j(t)).$

Since $\varphi_k(t_0, \cdot)$ maps $\mathbb{H} \setminus K_k(t_0)$ conformally onto $\mathbb{H}$, so from conformal invariance of extremal length, $(K_{j,t_0}(t), 0 \leq t < S_j)$ is also a Loewner chain in $\mathbb{H}$ w.r.t. $\infty$. Let $v_{j,t_0}(t) = \text{hcap}(K_{j,t_0}(t))/2$ for $0 \leq t < S_j$, and $L_{j,t_0}(t) = K_{j,t_0}(v_{j,t_0}^{-1}(t))$ for $0 \leq t < S_{j,t_0} := v_{j,t_0}(S_j)$. From Proposition 2.2, $L_{j,t_0}(t)$, $0 \leq t < S_{j,t_0}$, are chordal Loewner hulls driven by some $\eta_{j,t_0} \in C([0, S_{j,t_0}))$. Let $\psi_{j,t_0}(t,\cdot)$, $0 \leq t < S_{j,t_0}$, denote the corresponding chordal Loewner maps. Let $\xi_{j,t_0}(t) = \eta_{j,t_0}(v_{j,t_0}(t))$ and $\varphi_{j,t_0}(t,\cdot) = \psi_{j,t_0}(v_{j,t_0}(t),\cdot)$ for $0 \leq t < S_j$. Since $\psi_{j,t_0}(t,\cdot) = \varphi_{L_{j,t_0}(t)}$ for $0 \leq t < S_{j,t_0}$, so $\varphi_{j,t_0}(t,\cdot) = \varphi_{K_{j,t_0}(t)}$ for $0 \leq t < S_j$. We use $\partial_1$ and $\partial_z$ to denote the partial derivatives of $\varphi_j(\cdot,\cdot)$ and $\varphi_{j,t_0}(\cdot,\cdot)$ w.r.t. the first (real) and second (complex) variables, respectively, inside the bracket; and use $\partial_0$ to denote the partial derivative of $\varphi_{j,t_0}(\cdot,\cdot)$ w.r.t. the subscript $t_0$.

Fix $j \neq k \in \{1,2\}$, $t \in [0, S_j)$ and $s \in [0, S_k)$. Since $\varphi_k(s, \cdot) = \varphi_{K_k(s)}$, $\varphi_j(t, \cdot) = \varphi_{K_j(t)}$, $\varphi_{j,s}(t, \cdot) = \varphi_{K_{j,s}(t)}$ and $\varphi_{k,t}(s, \cdot) = \varphi_{K_{k,t}(s)}$, so from (3.1), for any $z \in \mathbb{H} \setminus (K_j(t) \cup K_k(s))$,

(3.2) $$\varphi_{K_j(t) \cup K_k(s)}(z) = \varphi_{k,t}(s, \varphi_j(t,z)) = \varphi_{j,s}(t, \varphi_k(s,z)).$$

Fix $\varepsilon \in (0, S_j - t)$. Since $K_{j,s}(r) = (K_j(r) \cup K_k(s))/K_k(s)$ for $r \in [0, S_j)$, so

$$\frac{L_{j,s}(v_{j,s}(t+\varepsilon))}{L_{j,s}(v_{j,s}(t))} = \frac{K_{j,s}(t+\varepsilon)}{K_{j,s}(t)} = \frac{K_j(t+\varepsilon) \cup K_k(s)}{K_j(t) \cup K_k(s)}$$

(3.3)
$$= \varphi_{K_j(t) \cup K_k(s)}(K_j(t+\varepsilon) \setminus K_j(t))$$
$$= \varphi_{k,t}(s, K_j(t+\varepsilon)/K_j(t)).$$

From Proposition 2.2 and (3.3), we have

(3.4) $$\{\xi_j(t)\} = \bigcap_{\varepsilon > 0} \overline{K_j(t+\varepsilon)/K_j(t)}$$

and

(3.5) $$\{\xi_{j,s}(t)\} = \{\eta_{j,s}(v_{j,s}(t))\} = \bigcap_{\varepsilon > 0} \overline{L_{j,s}(v_{j,s}(t+\varepsilon))/L_{j,s}(v_{j,s}(t))}$$

(3.6) $$= \bigcap_{\varepsilon > 0} \overline{(K_j(t+\varepsilon) \cup K_k(s))/(K_j(t) \cup K_k(s))}.$$

From (3.3)–(3.5), we have

(3.7) $$\xi_{j,s}(t) = \varphi_{k,t}(s, \xi_j(t)).$$

From Proposition 2.2 again, we have $\operatorname{hcap}(K_j(t+\varepsilon)/K_j(t)) = 2\varepsilon$ and

$$\operatorname{hcap}(L_{j,s}(v_{j,s}(t+\varepsilon))/L_{j,s}(v_{j,s}(t))) = 2(v_{j,s}(t+\varepsilon) - v_{j,s}(t)).$$

So from Proposition 2.1 and (3.3), we have

(3.8) $$v'_{j,s}(t) = \partial_z \varphi_{k,t}(s, \xi_j(t))^2.$$

Since $\varphi_{j,s}(t, z) = \psi_{j,s}(v_{j,s}(t), z)$, so for fixed $s \in [0, S_k)$, $(t, z) \mapsto \varphi_{j,s}(t, z)$ is $C^{1,a}$ differentiable, where the superscript "$a$" means analytic, and

(3.9)
$$\partial_1 \varphi_{j,s}(t, z) = \frac{2 v'_{j,s}(t)}{\psi_{j,s}(v_{j,s}(t), z) - \eta_{j,s}(v_{j,s}(t))}$$
$$= \frac{2 \partial_z \varphi_{k,t}(s, \xi_j(t))^2}{\varphi_{j,s}(t, z) - \varphi_{k,t}(s, \xi_j(t))}.$$

From (3.2), we see that $(s, t, z) \mapsto \varphi_{j,s}(t, z)$ is $C^{1,1,a}$ differentiable. Differentiate (3.9) using $\partial_z$, and then divide both sides by $\partial_z \varphi_{j,s}(t, z)$. We get

(3.10) $$\partial_1 \ln(\partial_z \varphi_{j,s}(t, z)) = \frac{-2 \partial_z \varphi_{k,t}(s, \xi_j(t))^2}{(\varphi_{j,s}(t, z) - \varphi_{k,t}(s, \xi_j(t)))^2}.$$

Differentiate (3.10) using $\partial_z$. We get

$$(3.11) \qquad \partial_1\left(\frac{\partial_z^2 \varphi_{j,s}(t,z)}{\partial_z \varphi_{j,s}(t,z)}\right) = \frac{4\,\partial_z \varphi_{k,t}(s,\xi_j(t))^2\, \partial_z \varphi_{j,s}(t,z)}{(\varphi_{j,s}(t,z) - \varphi_{k,t}(s,\xi_j(t)))^3}.$$

Differentiate (3.11) using $\partial_z$. We get

$$(3.12) \qquad \partial_1 \partial_z\left(\frac{\partial_z^2 \varphi_{j,s}(t,z)}{\partial_z \varphi_{j,s}(t,z)}\right) = \frac{4\,\partial_z \varphi_{k,t}(s,\xi_j(t))^2\, \partial_z^2 \varphi_{j,s}(t,z)}{(\varphi_{j,s}(t,z) - \varphi_{k,t}(s,\xi_j(t)))^3}$$
$$-\frac{12\,\partial_z \varphi_{k,t}(s,\xi_j(t))^2\, \partial_z \varphi_{j,s}(t,z)^2}{(\varphi_{j,s}(t,z) - \varphi_{k,t}(s,\xi_j(t)))^4}.$$

LEMMA 3.1. *For any $j \neq k \in \{0,1\}$, $t \in [0, S_j)$ and $s \in [0, S_k)$, we have*

$$(3.13) \quad \partial_0 \varphi_{k,t}(s, \xi_j(t)) = -3\,\partial_z^2 \varphi_{k,t}(s, \xi_j(t));$$

$$(3.14) \quad \frac{\partial_0\, \partial_z \varphi_{k,t}(s, \xi_j(t))}{\partial_z \varphi_{k,t}(s, \xi_j(t))} = \frac{1}{2} \cdot \left(\frac{\partial_z^2 \varphi_{k,t}(s, \xi_j(t))}{\partial_z \varphi_{k,t}(s, \xi_j(t))}\right)^2 - \frac{4}{3} \cdot \frac{\partial_z^3 \varphi_{k,t}(s, \xi_j(t))}{\partial_z \varphi_{k,t}(s, \xi_j(t))}.$$

PROOF. Differentiating both sides of the second "=" in (3.2) w.r.t. $t$, we get

$$\partial_0 \varphi_{k,t}(s, \varphi_j(t,z)) + \partial_z \varphi_{k,t}(s, \varphi_j(t,z))\, \partial_1 \varphi_j(t,z) = \partial_1 \varphi_{j,s}(t, \varphi_k(s,z))$$

for any $z \in \mathbb{H} \setminus (K_j(t) \cup K_k(s))$. So from (2.1), (3.2) and (3.9),

$$\partial_0 \varphi_{k,t}(s, \varphi_j(t,z)) = \frac{2\,\partial_z \varphi_{k,t}(s, \xi_j(t))^2}{\varphi_{k,t}(s, \varphi_j(t,z)) - \varphi_{k,t}(s, \xi_j(t))} - \frac{2\,\partial_z \varphi_{k,t}(s, \varphi_j(t,z))}{\varphi_j(t,z) - \xi_j(t)}$$

for any $z \in \mathbb{H} \setminus (K_j(t) \cup K_k(s))$. Since $\varphi_j(t, \cdot)$ maps $\mathbb{H} \setminus (K_j(t) \cup K_k(s))$ conformally onto $\mathbb{H} \setminus K_{k,t}(s)$, so for any $w \in \mathbb{H} \setminus K_{k,t}(s)$,

$$(3.15) \qquad \partial_0 \varphi_{k,t}(s, w) = \frac{2\,\partial_z \varphi_{k,t}(s, \xi_j(t))^2}{\varphi_{k,t}(s, w) - \varphi_{k,t}(s, \xi_j(t))} - \frac{2\,\partial_z \varphi_{k,t}(s, w)}{w - \xi_j(t)}.$$

In the above equation, let $w \to \xi_j(t)$ in $\mathbb{H} \setminus K_{k,t}(s)$. From the Taylor expansion of $\varphi_{k,t}(s, \cdot)$ at $\xi_j(t)$, we get (3.13). Differentiating (3.15) using $\partial_z$, we get

$$\partial_0\, \partial_z \varphi_{k,t}(s, w) = -\frac{2\,\partial_z \varphi_{k,t}(s, \xi_j(t))^2\, \partial_z \varphi_{k,t}(s, w)}{(\varphi_{k,t}(s, w) - \varphi_{k,t}(s, \xi_j(t)))^2}$$
$$-\frac{2\,\partial_z^2 \varphi_{k,t}(s, w)}{w - \xi_j(t)} + \frac{2\,\partial_z \varphi_{k,t}(s, w)}{(w - \xi_j(t))^2}.$$

Let $w \to \xi_j(t)$ in $\mathbb{H} \setminus K_{k,t}(s)$, then we get (3.14) from the Taylor expansion. $\square$

**4. Two-dimensional continuous local martingale.** Let $\kappa \in (0,4]$ and $x_1 < x_2 \in \mathbb{R}$. Let $X_1(t)$ and $X_2(t)$ be two independent Bessel process of dimension $3 - 8/\kappa$ started from $(x_2 - x_1)/\sqrt{\kappa}$. Let $T_j$ denote the first time that $X_j(t)$ visits 0, which exists and is finite because $3 - 8/\kappa \leq 1$. For $j = 1, 2$, let $Y_j(t) = \sqrt{\kappa} X_j(t)$, $0 \leq t \leq T_j$. Then there are two independent Brownian motions $B_1(t)$ and $B_2(t)$ such that, for $j = 1, 2$ and $0 \leq t \leq T_j$,

$$Y_j(t) = (x_2 - x_1) + (-1)^j \sqrt{\kappa} B_j(t) + \int_0^t \frac{\kappa - 4}{Y_j(s)} \, ds.$$

Fix $j \neq k \in \{1, 2\}$. For $0 \leq t \leq T_j$, let

$$\xi_j(t) = x_j + \sqrt{\kappa} B_j(t) + (-1)^j \int_0^t \frac{\kappa - 6}{Y_j(s)} \, ds,$$

$$p_j(t) = x_k - (-1)^j \int_0^t \frac{2}{Y_j(s)} \, ds.$$

Then $\xi_j(0) = x_j$, $p_j(0) = x_k$ and $\xi_j(t) - p_j(t) = (-1)^j Y_j(t)$, $0 \leq t \leq T_j$. Thus,

$$(4.1) \quad d\xi_j(t) = \sqrt{\kappa} \, dB_j(t) + \frac{\kappa - 6}{\xi_j(t) - p_j(t)} \, dt \quad \text{and} \quad dp_j(t) = \frac{2 \, dt}{p_j(t) - \xi_j(t)}$$

for $0 \leq t < T$. Let $K_j(t)$ and $\varphi_j(t, \cdot)$, $0 \leq t \leq T_j$, denote the chordal Loewner hulls and maps driven by $\xi_j(t)$, $0 \leq t \leq T_j$. Then $(K_j(t), 0 \leq t < T_j)$ are an SLE$(\kappa, \kappa - 6)$ process [5] started from $x_j$ with force point at $x_k$; $T_j$ is the first time that $x_k$ is swallowed by $K_j(t)$; and $\varphi_j(t, x_k) = p_j(t)$, $0 \leq t < T_j$. It is well known (e.g., [3]) that after a time-change, $(K_j(t), 0 \leq t < T_j)$ has the same distribution as a chordal SLE$_\kappa(\mathbb{H}; x_j \to x_k)$. Since $\kappa \leq 4$, so there is a crosscut $\beta_j(t)$, $0 \leq t \leq T_j$, in $\mathbb{H}$ from $x_j$ to $x_k$, such that $K_j(t) = \beta_j((0, t])$ for $0 \leq t < T_j$ [9]. Here a crosscut in $\mathbb{H}$ from $a \in \mathbb{R}$ to $b \in \mathbb{R}$ is a simple curve $\beta(t)$, $0 \leq t \leq T$, that satisfies $\beta(0) = a$, $\beta(T) = b$, and $\beta(t) \in \mathbb{H}$ for $0 < t < T$.

For $j = 1, 2$, let $(\mathcal{F}_t^j)$ denote the filtration generated by $(B_j(t))$. Then $(\xi_j)$ is $(\mathcal{F}_t^j)$-adapted, and $T_j$ is an $(\mathcal{F}_t^j)$-stopping time. Let

$$\mathcal{D} = \{(t_1, t_2) \in [0, T_1) \times [0, T_2) : \overline{K_1(t_1)} \cap \overline{K_2(t_2)} = \varnothing\}.$$

For $0 \leq t_k < T_k$, let $T_j(t_k) \in (0, T_j]$ be the maximal such that $\overline{K_j(t)} \cap \overline{K_k(t_k)} \neq \varnothing$ for $0 \leq t < T_j(t_k)$. Now we use the notation in the last section. Let $(t_1, t_2) \in \mathcal{D}$. Since $\varphi_{k,t_j}(t_k, \cdot) = \varphi_{K_{k,t_j}(t_k)}$, so $\varphi_{k,t_j}(t_k, \cdot)$ maps $\mathbb{H} \setminus K_{k,t_j}(t_k)$ conformally onto $\mathbb{H}$. By the Schwarz reflection principle, $\varphi_{k,t_j}(t_k, \cdot)$ extends conformally to $\Sigma_{K_{k,t_j}(t_k)}$, where for a hull $H$ in $\mathbb{H}$ w.r.t. $\infty$, $\Sigma_H = \mathbb{C} \setminus (H \cup \{\overline{z} : z \in H\} \cup [\inf(\overline{H} \cap \mathbb{R}), \sup(\overline{H} \cap \mathbb{R})])$ (cf. [17]). For $j \neq k \in \{0, 1\}$ and $h \in \mathbb{Z}_{\geq 0}$, let $A_{j,h}(t_1, t_2) = \partial_z^h \varphi_{k,t_j}(t_k, \xi_j(t_j))$. The definition makes sense since $\xi_j(t_j) \in \Sigma_{K_{k,t_j}(t_k)}$. Moreover, we have $A_{j,h} \in \mathbb{R}$ for any $h \geq 0$ since $\varphi_{k,t_j}(t_k, \cdot)$ is real valued on a real open interval containing $\xi_j(t_j)$. From (3.2), we see that $A_{j,0}(t_1, t_2) = \varphi_{K_1(t_1) \cup K_2(t_2)}(\beta_j(t_j))$, $j = 1, 2$. Since $K_1(t_1)$ lies to the left of $K_2(t_2)$, so $A_{1,0}(t_1, t_2) < A_{2,0}(t_1, t_2)$. Since $\varphi_{k,t_j}(t_k, \cdot)$ maps a part of the upper half plane to the upper half plane, so $A_{j,1}(t_1, t_2) > 0$, $j = 1, 2$. For $(t_1, t_2) \in \mathcal{D}$, define $E(t_1, t_2) = A_{2,0}(t_1, t_2) - A_{1,0}(t_1, t_2) > 0$,

$$(4.2) \quad N(t_1, t_2) = \frac{A_{1,1}(t_1, t_2) A_{2,1}(t_1, t_2)}{E(t_1, t_2)^2} = \frac{A_{1,1}(t_1, t_2) A_{2,1}(t_1, t_2)}{(A_{2,0}(t_1, t_2) - A_{1,0}(t_1, t_2))^2} > 0$$

and

$$(4.3) \quad M(t_1, t_2) = \left(\frac{N(t_1, t_2) N(0, 0)}{N(t_1, 0) N(0, t_2)}\right)^\alpha \exp\left(-\lambda \int_0^{t_1} \int_0^{t_2} 2N(s_1, s_2)^2 \, ds_2 \, ds_1\right) > 0,$$

where

$$(4.4) \quad \alpha = \alpha(\kappa) = \frac{6 - \kappa}{2\kappa}, \qquad \lambda = \lambda(\kappa) = \frac{(8 - 3\kappa)(6 - \kappa)}{2\kappa}.$$

Note that $M(t_1, 0) = M(0, t_2) = 1$ for any $0 \le t_1 < T_1$ and $0 \le t_2 < T_2$.

REMARK. If $\kappa < 8/3$, that is, $\lambda > 0$, then

$$\exp\left(-\lambda \int_0^{t_1} \int_0^{t_2} 2N(s_1, s_2)^2 \, ds_2 \, ds_1\right)$$

is the probability that in a loop soup [7] in $\mathbb{H}$ with intensity $\lambda$, there is no loop that intersects both $K_1(t_1)$ and $K_2(t_2)$.

THEOREM 4.1. (i) *For any fixed $(\mathcal{F}_t^2)$-stopping time $\bar{t}_2$ with $\bar{t}_2 < T_2$, $(M(t_1, \bar{t}_2), 0 \le t_1 < T_1(\bar{t}_2))$ is a continuous $(\mathcal{F}_{t_1}^1 \times \mathcal{F}_{\bar{t}_2}^2)_{t_1 \ge 0}$-local martingale, and*

$$(4.5) \quad \left.\frac{\partial_1 M}{M}\right|_{(t_1, \bar{t}_2)} = \left(3 - \frac{\kappa}{2}\right)\left(\left(\frac{A_{1,2}}{A_{1,1}} + \frac{2 A_{1,1}}{A_{2,0} - A_{1,0}}\right)\bigg|_{(t_1, \bar{t}_2)} - \frac{2}{p_1(t_1) - \xi_1(t_1)}\right)$$
$$\times \frac{\partial B_1(t_1)}{\sqrt{\kappa}}.$$

(ii) *For any fixed $(\mathcal{F}_t^1)$-stopping time $\bar{t}_1$ with $\bar{t}_1 < T_1$, $(M(\bar{t}_1, t_2), 0 \le t_2 < T_2(\bar{t}_1))$ is a continuous $(\mathcal{F}_{\bar{t}_1}^1 \times \mathcal{F}_{t_2}^2)_{t_2 \ge 0}$-local martingale, and*

$$(4.6) \quad \left.\frac{\partial_2 M}{M}\right|_{(\bar{t}_1, t_2)} = \left(3 - \frac{\kappa}{2}\right)\left(\left(\frac{A_{2,2}}{A_{2,1}} + \frac{2 A_{2,1}}{A_{1,0} - A_{2,0}}\right)\bigg|_{(\bar{t}_1, t_2)} - \frac{2}{p_2(t_2) - \xi_2(t_2)}\right)$$
$$\times \frac{\partial B_2(t_2)}{\sqrt{\kappa}}.$$

PROOF. Since $\varphi_{2,t_1}(0,\cdot) = \operatorname{id}_{\mathbb{H}}$, $\varphi_{1,0}(t_1,\cdot) = \varphi_1(t_1,\cdot)$, and $\xi_2(0) = x_2$, so

$$A_{1,0}(t_1,0) = \varphi_{2,t_1}(0,\xi_1(t_1)) = \xi_1(t_1), \qquad A_{1,1}(t_1,0) = 1;$$

$$A_{2,0}(t_1,0) = \varphi_{1,0}(t_1,\xi_2(0)) = \varphi_1(t_1,x_2) = p_1(t_1), \qquad A_{2,1}(t_1,0) = \partial_z\varphi_1(t_1,x_2).$$

Thus, $N(t_1,0) = \partial_z\varphi_1(t_1,x_2)/(p_1(t_1) - \xi_1(t_1))^2$. From the chordal Loewner equation, we get

$$\partial_{t_1}\partial_z\varphi_1(t_1,x_2) = \frac{-2\,\partial_z\varphi_1(t_1,x_2)}{(\varphi_1(t_1,x_2) - \xi_1(t_1))^2} = \frac{-2\,\partial_z\varphi_1(t_1,x_2)}{(p_1(t_1) - \xi_1(t_1))^2}.$$

From (4.1), we get

$$\partial_{t_1}(p_1(t_1) - \xi_1(t_1)) = -\partial\xi_1(t_1) + \frac{2\,\partial t_1}{p_1(t_1) - \xi_1(t_1)}.$$

From the above two formulas and Itô's formula, we get

(4.7) $\quad \partial_1 N(t_1,0)^\alpha / (\alpha N(t_1,0)^\alpha) = 2\,\partial\xi_1(t_1)/(p_1(t_1) - \xi_1(t_1)).$

Now fix an $(\mathcal{F}_t^2)$-stopping time $\bar{t}_2$ with $\bar{t}_2 < T_2$. Then we get a filtration $(\mathcal{F}_t^1 \times \mathcal{F}_{\bar{t}_2}^2)_{t\geq 0}$. Since $B_1(t)$ and $B_2(t)$ are independent, so $B_1(t)$ is an $(\mathcal{F}_t^1 \times \mathcal{F}_{\bar{t}_2}^2)_{t\geq 0}$-Brownian motion. Then $T_1(\bar{t}_2)$ is an $(\mathcal{F}_t^1 \times \mathcal{F}_{\bar{t}_2}^2)_{t\geq 0}$-stopping time, $A_{j,h}(t,\bar{t}_2)$, $j=1,2$, $E(t,\bar{t}_2)$, $N(t,\bar{t}_2)$ and $M(t,\bar{t}_2)$ are defined for $t \in [0, T_1(\bar{t}_2))$. From the chordal Loewner equation and (3.2), $\varphi_{1,\bar{t}_2}(t,\cdot)$ and $\varphi_{2,t}(\bar{t}_2,\cdot)$, $0 \leq t < T_1(\bar{t}_2)$, are $(\mathcal{F}_t^1 \times \mathcal{F}_{\bar{t}_2}^2)_{t\geq 0}$-adapted. Since $A_{1,h}(t,\bar{t}_2) = \partial_z^h \varphi_{2,t}(\bar{t}_2, \xi_1(t))$, so from Itô's formula, $(A_{1,h}(t_1,\bar{t}_2), 0 \leq t_1 < T_1(\bar{t}_2))$ satisfies the $(\mathcal{F}_t^1 \times \mathcal{F}_{\bar{t}_2}^2)_{t\geq 0}$-adapted SDE:

(4.8)
$$\partial_1 A_{1,h}(t_1,\bar{t}_2) = A_{1,h+1}(t_1,\bar{t}_2)\,\partial\xi_1(t_1)$$
$$+ \left(\partial_0\,\partial_z^h \varphi_{2,t_1}(\bar{t}_2,\xi_1(t_1)) + \frac{\kappa}{2} A_{1,h+2}(t,\bar{t}_2)\right)\partial t_1.$$

From (3.9) and (3.10), we have

(4.9)
$$\partial_1 A_{2,0}(t_1,t_2) = \frac{2A_{1,1}(t_1,t_2)^2}{E(t_1,t_2)}\,\partial t_1,$$
$$\frac{\partial_1 A_{2,1}(t_1,t_2)}{A_{2,1}(t_1,t_2)} = -\frac{2A_{1,1}(t_1,t_2)^2}{E(t_1,t_2)^2}\,\partial t_1.$$

From (4.8), (4.9) and Lemma 3.1, we have

(4.10) $\quad \partial_1 A_{1,0} = A_{1,1}\,\partial\xi_1(t_1) + \left(\frac{\kappa}{2} - 3\right)A_{1,2}\,\partial t_1$

and

(4.11) $\quad \dfrac{\partial_1 A_{1,1}}{A_{1,1}} = \dfrac{A_{1,2}}{A_{1,1}}\,\partial\xi_1(t_1) + \left(\dfrac{1}{2}\cdot\left(\dfrac{A_{1,2}}{A_{1,1}}\right)^2 + \left(\dfrac{\kappa}{2} - \dfrac{4}{3}\right)\cdot\dfrac{A_{1,3}}{A_{1,1}}\right)\partial t_1,$

where "$(t_1, \bar{t}_2)$" are omitted. Since $E = A_{2,0} - A_{1,0}$, from (4.9) and (4.10), we have

$$(4.12) \qquad \partial_1 E = -A_{1,1} \partial \xi_1(t_1) + \left( \frac{2A_{1,1}^2}{E} + \left(3 - \frac{\kappa}{2}\right) A_{1,2} \right) \partial t_1.$$

Let $C_h = A_{1,h}$ for $h \in \mathbb{Z}_{\geq 0}$. From (4.9)–(4.12) and Itô's formula, we have

$$(4.13) \qquad \frac{\partial_1 N^\alpha}{\alpha N^\alpha} = \left( \frac{C_2}{C_1} + \frac{2C_1}{E} \right) \partial \xi_1(t_1) + (8 - 3\kappa)\left( \frac{1}{4} \cdot \frac{C_2^2}{C_1^2} - \frac{1}{6} \cdot \frac{C_3}{C_1} \right) \partial t_1.$$

The above SDE is $(\mathcal{F}_t^1 \times \mathcal{F}_{\bar{t}_2}^2)_{t \geq 0}$-adapted. Now (4.7) is also an $(\mathcal{F}_t^1 \times \mathcal{F}_{\bar{t}_2}^2)_{t \geq 0}$-adapted SDE since $B_1(t)$ is an $(\mathcal{F}_t^1 \times \mathcal{F}_{\bar{t}_2}^2)_{t \geq 0}$-Brownian motion. Thus, from (4.1), (4.7), (4.13) and Itô's formula, we have

$$(4.14) \quad \begin{aligned} &\frac{\partial_1 (N(t_1, \bar{t}_2)/N(t_1, 0))^\alpha}{\alpha (N(t_1, \bar{t}_2)/N(t_1, 0))^\alpha} \\ &= \left( \frac{C_2(t_1, \bar{t}_2)}{C_1(t_1, \bar{t}_2)} + \frac{2C_1(t_1, \bar{t}_2)}{E(t_1, \bar{t}_2)} - \frac{2}{p_1(t_1) - \xi_1(t_1)} \right) \sqrt{\kappa} \partial B_1(t_1) \\ &\quad + (8 - 3\kappa) \left( \frac{1}{4} \cdot \frac{C_2(t_1, \bar{t}_2)^2}{C_1(t_1, \bar{t}_2)^2} - \frac{1}{6} \cdot \frac{C_3(t_1, \bar{t}_2)}{C_1(t_1, \bar{t}_2)} \right) \partial t_1. \end{aligned}$$

Since $C_j(t_1, t_2) = \partial_z^j \varphi_{2,t_1}(t_2, \xi_1(t_1))$, so $\partial_2 C_j(t_1, t_2) = \partial_1 \partial_z^j \varphi_{2,t_1}(t_2, \xi_1(t_1))$, and

$$\left( \frac{1}{4} \cdot \frac{C_2^2}{C_1^2} - \frac{1}{6} \cdot \frac{C_3}{C_1} \right)\bigg|_{(t_1,t_2)} = \frac{1}{12}(\partial_z^2/\partial_z)\varphi_{2,t_1}(t_2, \xi_1(t_1))^2$$

$$- \frac{1}{6}\partial_z(\partial_z^2/\partial_z)\varphi_{2,t_1}(t_2, \xi_1(t_1)).$$

From (3.11) and (3.12), we have

$$\frac{\partial}{\partial t_2}[(\partial_z^2/\partial_z)\varphi_{2,t_1}(t_2, \xi_1(t_1))^2] = \frac{8A_{2,1}^2 C_2}{E^3}\bigg|_{(t_1,t_2)},$$

$$\frac{\partial}{\partial t_2}[\partial_z(\partial_z^2/\partial_z)\varphi_{2,t_1}(t_2, \xi_1(t_1))] = \left( \frac{4A_{2,1}^2 C_2}{E^3} - \frac{12A_{2,1}^2 C_1^2}{E^4} \right)\bigg|_{(t_1,t_2)}.$$

From the above three formulas, we get

$$\partial_2 \left( \frac{1}{4} \cdot \frac{C_2^2}{C_1^2} - \frac{1}{6} \cdot \frac{C_3}{C_1} \right)\bigg|_{(t_1,t_2)} = \frac{2A_{2,1}^2 C_1^2}{E^4}\bigg|_{(t_1,t_2)} = 2N(t_1, t_2)^2.$$

Since $\varphi_{2,t_1}(0, \cdot) = \mathrm{id}_\mathbb{H}$, so $\partial_z^j \varphi_{2,t_1}(0, \cdot) = 0$ for $j \geq 2$. Thus, $C_2(t_1, 0) = C_3(t_1, 0) = 0$. So

$$(4.15) \qquad \frac{1}{4} \cdot \frac{C_2(t_1, t_2)^2}{C_1(t_1, t_2)^2} - \frac{1}{6} \cdot \frac{C_3(t_1, t_2)}{C_1(t_1, t_2)} = \int_0^{t_2} 2N(t_1, s_2)^2 \, ds_2.$$

Then (4.5) follows from (4.3)–(4.4) and (4.14)–(4.15); (4.6) follows from the symmetry. □

Now we make some improvement over the above theorem. Let $\bar{t}_2$ be an $(\mathcal{F}_t^2)$-stopping time with $\bar{t}_2 < T_2$. Suppose $R$ is an $(\mathcal{F}_t^1 \times \mathcal{F}_{\bar{t}_2}^2)_{t \geq 0}$-stopping time with $R < T_1(\bar{t}_2)$. Let $\mathcal{F}_{R,\bar{t}_2}$ denote the $\sigma$-algebra obtained from the filtration $(\mathcal{F}_t^1 \times \mathcal{F}_{\bar{t}_2}^2)_{t \geq 0}$ and its stopping time $R$, that is, $\mathcal{E} \in \mathcal{F}_{R,\bar{t}_2}$ iff for any $t \geq 0$, $\mathcal{E} \cap \{R \leq t\} \in \mathcal{F}_t^1 \times \mathcal{F}_{\bar{t}_2}^2$. For every $t \geq 0$, $R + t$ is also an $(\mathcal{F}_t^1 \times \mathcal{F}_{\bar{t}_2}^2)_{t \geq 0}$-stopping time. So we have a filtration $(\mathcal{F}_{R+t,\bar{t}_2})_{t \geq 0}$. Since $(\xi_1(t))$ and $(p_1(t))$ are $(\mathcal{F}_t^1 \times \mathcal{F}_{\bar{t}_2}^2)_{t \geq 0}$-adapted, so $(\xi_1(R+t), t \geq 0)$, $(p_1(R+t), t \geq 0)$, $(\varphi_1(R+t, \cdot), t \geq 0)$ and $(K_1(R+t), t \geq 0)$ are $(\mathcal{F}_{R+t,\bar{t}_2})_{t \geq 0}$-adapted. Suppose $I \in [0, \bar{t}_2]$ is $\mathcal{F}_{R,\bar{t}_2}$-measurable. From $I \leq \bar{t}_2$ we have $T_1(I) \geq T_1(\bar{t}_2) > R$. Then $\varphi_{1,I}(R+t, \cdot)$ and $\varphi_{2,R+t}(I, \cdot)$ are defined for $0 \leq t < T_1(I) - R$.

LEMMA 4.1. $T_1(I) - R$ is an $(\mathcal{F}_{R+t,\bar{t}_2})_{t \geq 0}$-stopping time and $(\varphi_{1,I}(R+t, \cdot), 0 \leq t < T_1(I) - R)$ and $(\varphi_{2,R+t}(I, \cdot), 0 \leq t < T_1(I) - R)$ are $(\mathcal{F}_{R+t,\bar{t}_2})_{t \geq 0}$-adapted.

PROOF. Since $T_1(I) - R > t$ iff $K_1(R+t) \cap K_2(I) = \varnothing$, and that $(\varphi_1(R+t, \cdot))$, and $(K_1(R+t))$ are $\mathcal{F}_{R+t,\bar{t}_2}$-adapted, so from (3.2), we suffice to show that $\varphi_2(I, \cdot)$ is $\mathcal{F}_{R,\bar{t}_2}$-measurable. Fix $n \in \mathbb{N}$. Let $I_n = \lfloor nI \rfloor / n$. For $m \in \mathbb{N} \cup \{0\}$, let $\mathcal{E}_n(m) = \{m/n \leq I_n < (m+1)/n\}$. Then $\mathcal{E}_n(m)$ is $\mathcal{F}_{R,\bar{t}_2}$-measurable, and $I_n = m/n$ on $\mathcal{E}_n(m)$. Since $m/n \leq \bar{t}_2$ and $I_n = m/n$ on $\mathcal{E}_n(m)$, so $I_n$ agrees with $(m/n) \wedge \bar{t}_2$ on $\mathcal{E}_n(m)$. Now $(m/n) \wedge \bar{t}_2$ is an $(\mathcal{F}_t^2)$-stopping time, and $\mathcal{F}_{(m/n) \wedge \bar{t}_2}^2 \subset \mathcal{F}_{\bar{t}_2}^2 \subset \mathcal{F}_{R,\bar{t}_2}$. So $\varphi_2((m/n) \wedge \bar{t}_2, \cdot)$ is $\mathcal{F}_{R,\bar{t}_2}$-measurable. Since $\varphi_2(I_n, \cdot) = \varphi_2((m/n) \wedge \bar{t}_2, \cdot)$ on $\mathcal{E}_n(m)$, and $\mathcal{E}_n(m)$ is $\mathcal{F}_{R,\bar{t}_2}$-measurable for each $m \in \mathbb{N} \cup \{0\}$, so $\varphi_2(I_n, \cdot)$ is $\mathcal{F}_{R,\bar{t}_2}$-measurable. Since $\varphi_2(I_n, \cdot) \to \varphi_2(I, \cdot)$ as $n \to \infty$, so $\varphi_2(I, \cdot)$ is also $\mathcal{F}_{R,\bar{t}_2}$-measurable. Then we are done. □

Let $B_1^R(t) = B_1(R+t) - B_1(R)$, $0 \leq t < \infty$. Since $B_1(t)$ is an $(\mathcal{F}_t^1 \times \mathcal{F}_{\bar{t}_2}^2)_{t \geq 0}$-Brownian motion, so $B_1^R(t)$ is an $(\mathcal{F}_{R+t,\bar{t}_2})_{t \geq 0}$-Brownian motion. Then $(\xi_1(R+t))$ satisfies the $(\mathcal{F}_{R+t,\bar{t}_2})_{t \geq 0}$-adapted SDE:

$$d\xi_1(R+t) = \sqrt{\kappa}\, dB_1^R(t) + \frac{\kappa - 6}{\xi_1(R+t) - p_1(R+t)}\, dt.$$

The SDEs in the proof of Theorem 4.1 still hold if $t_1$ is replaced by $R+t$, $\bar{t}_2$ is replaced by $I$, and $B_1(t_1)$ is replaced by $B_1^R(t_1)$. The difference is that the SDEs now are all $(\mathcal{F}_{R+t,\bar{t}_2})_{t \geq 0}$-adapted. So we have the following theorem.

THEOREM 4.2. (i) Suppose $\bar{t}_2$ is an $(\mathcal{F}_t^2)$-stopping time with $\bar{t}_2 < T_2$. Suppose $R$ is an $(\mathcal{F}_t^1 \times \mathcal{F}_{\bar{t}_2}^2)_{t \geq 0}$-stopping time with $R < T_1(\bar{t}_2)$. Let $I \in [0, \bar{t}_2]$ be $\mathcal{F}_{R,\bar{t}_2}$-measurable. Then $(M(R+t,I), 0 \leq t < T_1(I) - R)$ is a continuous $(\mathcal{F}_{R+t,\bar{t}_2})_{t\geq 0}$-local martingale.

(ii) Suppose $\bar{t}_1$ is an $(\mathcal{F}_t^1)$-stopping time with $\bar{t}_1 < T_1$. Suppose $I$ is an $(\mathcal{F}_{\bar{t}_1}^1 \times \mathcal{F}_t^2)_{t\geq 0}$-stopping time with $I < T_2(\bar{t}_1)$. Let $R \in [0, \bar{t}_1]$ be $\mathcal{F}_{\bar{t}_1, I}$-measurable. Then $(M(R, I+t), 0 \leq t < T_2(R) - I)$ is a continuous $(\mathcal{F}_{\bar{t}_1, I+t})_{t\geq 0}$-local martingale.

PROOF. (i) follows from the above argument. (ii) follows from the symmetry. □

**5. Boundedness.** We now use the notation and results in Section 5.2 of [17]. Let $H$ be a nonempty hull in $\mathbb{H}$ w.r.t. $\infty$. Then $a_H = \inf\{\overline{H} \cap \mathbb{R}\}$, $b_H = \sup\{\overline{H} \cap \mathbb{R}\}$, $\Sigma_H = \mathbb{C} \setminus (H \cup \{\overline{z} : z \in H\} \cup [a_H, b_H])$, and $\mathcal{H}(H)$ is the set of hulls in $\mathbb{H}$ w.r.t. $\infty$ that are contained in $H$. From Lemma 5.4 in [17], any sequence $(K_n)$ in $\mathcal{H}(H)$ contains a subsequence $(L_n)$ such that $\varphi_{L_n} \xrightarrow{\text{l.u.}} \varphi_K$ (converges locally uniformly) in $\Sigma_H$ for some $K \in \mathcal{H}(H)$. We now make some improvement over this result. Let $Q_H = \overline{H} \cap \mathbb{R}$. Then $Q_H$ is a closed subset of $[a_H, b_H]$. Let

$$\Sigma_H^* = \Sigma_H \cup ([a_H, b_H] \setminus Q_H) = \mathbb{C} \setminus (H \cup \{\overline{z} : z \in H\} \cup Q_H),$$

which may strictly contain $\Sigma_H$. For any $K \in \mathcal{H}(H)$, $\varphi_K$ extends conformally to $\Sigma_H^*$ by the Schwarz reflection principle, and $\varphi_K'(x) > 0$ for any $x \in \mathbb{R} \setminus Q_H$ from (5.1) in [17], so $\varphi_K$ preserves the order on $\mathbb{R} \setminus Q_H$.

LEMMA 5.1. Suppose $(K_n)$ is a sequence in $\mathcal{H}(H)$. Then it contains some subsequence $(L_n)$ such that $\varphi_{L_n} \xrightarrow{\text{l.u.}} \varphi_K$ in $\Sigma_H^*$ for some $K \in \mathcal{H}(H)$.

PROOF. From the argument after Corollary 5.1 in [17], there is $M_H > 0$ such that $|\varphi_K(z) - z| \leq M_H$ for any $K \in \mathcal{H}(H)$ and $z \in \Sigma_H$. After the extension, we have $|\varphi_K(z) - z| \leq M_H$ for any $K \in \mathcal{H}(H)$ and $z \in \Sigma_H^*$. So $\{\varphi_{K_n}(z) - z : n \in \mathbb{N}\}$ is a normal family in $\Sigma_H^*$. Then $(K_n)$ contains a subsequence $(L_n)$ such that $\varphi_{L_n}(z) - z \xrightarrow{\text{l.u.}} f(z)$ in $\Sigma_H^*$ for some $f$ that is analytic in $\Sigma_H^*$. So $\varphi_{L_n} \xrightarrow{\text{l.u.}} g$ in $\Sigma_H^*$, where $g(z) := z + f(z)$ is analytic in $\Sigma_H^*$. From Lemma 5.4 in [17], we may assume that $\varphi_{L_n} \xrightarrow{\text{l.u.}} \varphi_K$ in $\Sigma_H$ for some $K \in \mathcal{H}(H)$. Thus, $g = \varphi_K$ in $\Sigma_H$. Since they are both analytic in $\Sigma_H^*$, so $g = \varphi_K$ in $\Sigma_H^*$. Thus, $\varphi_{L_n} \xrightarrow{\text{l.u.}} \varphi_K$ in $\Sigma_H^*$. □

LEMMA 5.2. If $y_1 < y_2 < a_H$ or $y_1 > y_2 > b_H$, then $\varphi_H'(y_1) > \varphi_H'(y_2)$.

PROOF. This follows from differentiating (5.1) in [17] for $z \in \mathbb{R} \setminus [c_H, d_H]$, and the fact that $\varphi_H$ is increasing on $(-\infty, a_H)$ and $(b_H, \infty)$, and maps them to $(-\infty, c_H)$ and $(d_H, \infty)$, respectively. □

Let HP denote the set of $(H_1, H_2)$ such that $H_j$ is a hull in $\mathbb{H}$ w.r.t. $\infty$ that contains some neighborhood of $x_j$ in $\mathbb{H}$, $j = 1, 2$, and $\overline{H_1} \cap \overline{H_2} = \varnothing$. Let $(H_1, H_2) \in \text{HP}$. Then $b_{H_1} < a_{H_2}$, $H_1 \cup H_2$ is a hull in $\mathbb{H}$ w.r.t. $\infty$, and $Q_{H_1 \cup H_2} = Q_{H_1} \cup Q_{H_2} \subset [a_{H_1}, b_{H_1}] \cup [a_{H_2}, b_{H_2}]$. Let $T_j(H_j)$ be the first time that $\overline{K_j(t)} \cap \overline{\mathbb{H} \setminus H_j} \ne \varnothing$, $j = 1, 2$. Then $T_j(H_j)$ is an $(\mathcal{F}_t^j)$-stopping time, $0 < T_j(H_j) < T_j$, and $K_j(t) \subset H_j$ for $0 \le t \le T_j(H_j)$. Thus,

$$(5.1) \qquad T_j(H_j) = \text{hcap}(K_j(T_j(H_j)))/2 \le \text{hcap}(H_j)/2.$$

THEOREM 5.1. *For any $(H_1, H_2) \in \text{HP}$, there are $C_2 > C_1 > 0$ depending only on $H_1$ and $H_2$ such that $M(t_1, t_2) \in [C_1, C_2]$ for any $(t_1, t_2) \in [0, T_1(H_1)] \times [0, T_2(H_2)]$.*

PROOF. Let $(H_1, H_2) \in \text{HP}$ and $H = H_1 \cup H_2$. Throughout this proof, we use $C_n$, $n \in \mathbb{N}$, to denote some positive constant that depends only on $H_1$ and $H_2$. From (4.3) and (5.1), we suffice to show that for some $C_4 > C_3 > 0$, $N(t_1, t_2) \in [C_3, C_4]$ for $(t_1, t_2) \in [0, T_1(H_1)] \times [0, T_2(H_2)]$. Fix $(t_1, t_2) \in [0, T_1(H_1)] \times [0, T_2(H_2)]$. First suppose $t_1, t_2 > 0$. Fix $j \ne k \in \{1, 2\}$. For any $s_j \in [0, t_j)$, from (3.4) we have $\xi_j(s_j) \in \overline{K_j(t_j)/K_j(s_j)}$, so

$$\xi_j(s_j) \in [a_{K_j(t_j)/K_j(s_j)}, b_{K_j(t_j)/K_j(s_j)}]$$
$$\subset [c_{K_j(t_j)/K_j(s_j)}, d_{K_j(t_j)/K_j(s_j)}] \subset [c_{K_j(t_j)}, d_{K_j(t_j)}],$$

where the second and third inclusions follow from Lemma 5.2 and Lemma 5.3 in [17]. Let $s_j \to t_j$. We get $\xi_j(t_j) \in [c_{K_j(t_j)}, d_{K_j(t_j)}]$. For $s_j \in [0, t_j)$, from (3.6) and (3.7),

$$A_{j,0}(s_j, t_k) = \varphi_{k,s_j}(t_k, \xi_j(s_j)) \in \overline{(K_j(t_j) \cup K_k(t_k))/(K_j(s_j) \cup K_k(t_k))},$$

which implies that $A_{j,0}(s_j, t_k) \in [c_{K_j(t_j) \cup K_k(t_k)}, d_{K_j(t_j) \cup K_k(t_k)}] \subset [c_H, d_H]$. Let $s_j \to t_j$. We get $A_{j,0}(t_j, t_k) \in [c_H, d_H]$. This also holds for $A_{k,0}(t_j, t_k)$. Thus,

$$(5.2) \qquad |E(t_j, t_k)| = |A_{j,0}(t_j, t_k) - A_{k,0}(t_j, t_k)| \le d_H - c_H.$$

Fix $q_1, q_2, r_1, r_2 \in \mathbb{R}$ with $r_1 < a_{H_1} \le b_{H_1} < q_1 < q_2 < a_{H_2} \le b_{H_2} < r_2$. From Lemma 5.1, there are $C_6 > C_5 > 0$ such that, for $x = q_1, q_2, r_1, r_2$, $\partial_z \varphi_{K_1(t_1) \cup K_2(t_2)}(x)$, $\partial_z \varphi_1(t_1, x)$ and $\partial_z \varphi_2(t_2, x)$ all lie in $[C_5, C_6]$. Fix $j \ne k \in \{1, 2\}$. From (3.2) there are $C_8 > C_7 > 0$ such that, for $x = q_j, r_j$, $\partial_z \varphi_{k,t_j}(t_k, \varphi_j(t_j, x)) \in [C_7, C_8]$. Since $[a_{K_j(t_j)}, b_{K_j(t_j)}] \subset [a_{H_j}, b_{H_j}]$, so $r_j$ is disconnected from $q_j$ in $\mathbb{R}$ by $[a_{K_j(t_j)}, b_{K_j(t_j)}]$. Since $\varphi_j(t_j, \cdot) = \varphi_{K_j(t_j)}$, so $\varphi_j(t_j, r_j)$ is disconnected from $\varphi_j(t_j, q_j)$ in $\mathbb{R}$ by $[c_{K_j(t_j)}, d_{K_j(t_j)}]$. Since $\xi_j(t_j) \in [c_{K_j(t_j)}, d_{K_j(t_j)}]$, so $\xi_j(t_j)$ lies between $\varphi_j(t_j, r_j)$ and $\varphi_j(t_j, q_j)$. Since $r_j$ and $q_j$ lie on the same side of $K_k(t_k)$, so $\varphi_j(t_j, r_j)$, $\xi_j(t_j)$, and $\varphi_j(t_j, q_j)$ lie on the same side of $\varphi_j(t_j, K_k(t_k)) = K_{k,t_j}(t_k)$. Since $\varphi_{k,t_j}(t_k, \cdot) = \varphi_{K_{k,t_j}(t_k)}$, so from Lemma 5.2, $\partial_z \varphi_{k,t_j}(t_k, \xi_j(t_j))$ lies between $\partial_z \varphi_{k,t_j}(t_k, \varphi_j(t_j, r_j))$ and $\partial_z \varphi_{k,t_j}(t_k, \varphi_j(t_j, q_j))$. Thus,

(5.3) $$A_{j,1}(t_j, t_k) = \partial_z \varphi_{k,t_j}(t_k, \xi_j(t_j)) \in [C_7, C_8].$$

From (3.2) and the above argument, we see that $A_{j,0}(t_j, t_k) = \varphi_{k,t_j}(t_k, \xi_j(t_j))$ lies between $\varphi_{K_j(t_j) \cup K_k(t_k)}(r_j)$ and $\varphi_{K_j(t_j) \cup K_k(t_k)}(q_j)$ for $j = 1, 2$. Since $r_1 < q_1 < q_2 < r_2$, so

$$\varphi_{K_1(t_2) \cup K_2(t_2)}(r_1) < \varphi_{K_1(t_1) \cup K_2(t_2)}(q_1) < \varphi_{K_1(t_1) \cup K_2(t_2)}(q_2) < \varphi_{K_1(t_1) \cup K_2(t_2)}(r_2).$$

From Lemma 5.1, there is $C_9 > 0$ such that $\partial_z \varphi_{K_1(t_1) \cup K_2(t_2)}(x) \geq C_9$ for $x \in [q_1, q_2]$. So

(5.4) $$|E(t_1, t_2)| \geq \varphi_{K_1(t_1) \cup K_2(t_2)}(q_2) - \varphi_{K_1(t_1) \cup K_2(t_2)}(q_1) \geq C_9(q_2 - q_1).$$

From (5.2), (5.3) and (5.4), we have $C_4 > C_3 > 0$ such that $N(t_1, t_2) \in [C_3, C_4]$ for $(t_1, t_2) \in (0, T_1(H_1)] \times (0, T_2(H_2)]$. By letting $t_1$ or $t_2$ tend to 0, we obtain the above inequality in the case $t_1$ or $t_2$ equals to 0. So we are done. □

Now we explain the meaning of $M(t_1, t_2)$. Fix $(H_1, H_2) \in \text{HP}$. Let $\mu$ denote the joint distribution of $(\xi_1(t) : 0 \leq t \leq T_1)$ and $(\xi_2(t) : 0 \leq t \leq T_2)$. From Theorem 4.1 and Theorem 5.1, we have $\int M(T_1(H_1), T_2(H_2)) d\mu = \mathbf{E}[M(T_1(H_1), T_2(H_2))] = M(0, 0) = 1$. Note that $M(T_1(H_1), T_2(H_2)) > 0$. Suppose $\nu$ is a measure on $\mathcal{F}^1_{T_1(H_1)} \times \mathcal{F}^2_{T_2(H_2)}$ such that $d\nu/d\mu = M(T_1(H_1), T_2(H_2))$. Then $\nu$ is a probability measure. Now suppose the joint distribution of $(\xi_1(t), 0 \leq t \leq T_1(H_1))$ and $(\xi_2(t), 0 \leq t \leq T_2(H_2))$ is $\nu$ instead of $\mu$. Fix an $(\mathcal{F}^2_t)$-stopping time $\bar{t}_2$ with $\bar{t}_2 \leq T_2(H_2)$. From (4.1), (4.5) and the Girsanov theorem [8], there is an $(\mathcal{F}^1_t \times \mathcal{F}^2_{\bar{t}_2})$-Brownian motion $\widetilde{B}_1(t)$ such that $\xi_1(t_1)$ satisfies the $(\mathcal{F}^1_{t_1} \times \mathcal{F}^2_{\bar{t}_2})$-adapted SDE for $0 \leq t_1 \leq T_1(H_1)$:

(5.5) $$d\xi_1(t_1) = \sqrt{\kappa} \, d\widetilde{B}_1(t_1)$$
$$+ \left(3 - \frac{\kappa}{2}\right) \left(\frac{A_{1,2}(t_1, \bar{t}_2)}{A_{1,1}(t_1, \bar{t}_2)} + \frac{2 A_{1,1}(t_1, \bar{t}_2)}{A_{2,0}(t_1, \bar{t}_2) - A_{1,0}(t_1, \bar{t}_2)}\right) dt_1.$$

From (4.10) and (5.5), we have

(5.6) $$dA_{1,0}(t_1, \bar{t}_2) = A_{1,1}(t_1, \bar{t}_2) \sqrt{\kappa} \, d\widetilde{B}_1(t) + \frac{(6 - \kappa) A_{1,1}(t_1, \bar{t}_2)^2 \, dt_1}{A_{2,0}(t_1, \bar{t}_2) - A_{1,0}(t_1, \bar{t}_2)}.$$

Recall that $A_{1,0}(t_1, \bar{t}_2) = \varphi_{2,t_1}(\bar{t}_2, \xi_1(t_1)) = \xi_{1,\bar{t}_2}(t_1) = \eta_{1,\bar{t}_2}(v_{1,\bar{t}_2}(t_1))$, and $v'_{1,\bar{t}_2}(t_1) = A_{1,1}(t_1, \bar{t}_2)^2$ [see (3.8)]. From (5.6), there is a Brownian motion $\widehat{B}_1(t_1)$ such that

(5.7) $$d\eta_{1,\bar{t}_2}(s_1) = \sqrt{\kappa} \, d\widehat{B}_1(s_1) + \frac{(\kappa - 6) \, ds_1}{\eta_{1,\bar{t}_2}(s_1) - A_{2,0}(v^{-1}_{1,\bar{t}_2}(s_1), \bar{t}_2)}.$$

Since $A_{2,0}(v_{1,\bar{t}_2}^{-1}(s_1), \bar{t}_2) = \varphi_{1,\bar{t}_2}(v_{1,\bar{t}_2}^{-1}(s_1), \xi_2(\bar{t}_2)) = \psi_{1,\bar{t}_2}(s_1, \xi_2(\bar{t}_2))$ and $\psi_{1,\bar{t}_2}(s, \cdot)$, $0 \leq s \leq v_{1,\bar{t}_2}(T_1(H_1))$, are chordal Loewner maps driven by $\eta_{1,\bar{t}_2}(s)$, so the chordal Loewner hulls $L_{1,\bar{t}_2}(s)$, $0 \leq s \leq v_{1,\bar{t}_2}(T_1(H_1))$, driven by $\eta_{1,\bar{t}_2}(s)$ are a part of the chordal SLE$(\kappa, \kappa - 6)$ process started from $\eta_{1,\bar{t}_2}(0) = \varphi_2(\bar{t}_2, x_1)$ with force point at $A_{2,0}(v_{1,\bar{t}_2}^{-1}(0), \bar{t}_2) = \xi_2(\bar{t}_2)$. Thus, after a time-change, it is a chordal SLE$_\kappa$ in $\mathbb{H}$ from $\varphi_2(\bar{t}_2, x_1)$ to $\xi_2(\bar{t}_2)$. Note that $\varphi_2(\bar{t}_2, \cdot)^{-1}$ maps $\mathbb{H}$ conformally onto $\mathbb{H} \setminus \beta_2((0, \bar{t}_2])$, maps $L_{1,\bar{t}_2}(v_{1,\bar{t}_2}(t_1))$ onto $K_1(t_1) = \beta_1((0, t_1])$, and takes $\varphi_2(\bar{t}_2, x_1)$ and $\xi_2(\bar{t}_2)$ to $x_1$ and $\beta_2(\bar{t}_2)$, respectively. Thus, $\beta_1(t)$, $0 \leq t \leq T_1(H_1)$, is the time-change of a chordal SLE$_\kappa$ trace in $\mathbb{H} \setminus \beta_2((0, \bar{t}_2])$ from $x_1$ to $\beta_2(\bar{t}_2)$, stopped on hitting $\overline{\mathbb{H} \setminus H_1}$. Similarly, for any $(\mathcal{F}_t^1)$-stopping time $\bar{t}_1$ with $\bar{t}_1 \leq T_1(H_1)$, $\beta_2(t)$, $0 \leq t \leq T_2(H_2)$, is a time-change of a chordal SLE$_\kappa$ trace in $\mathbb{H} \setminus \beta_1((0, \bar{t}_1])$ from $x_2$ to $\beta_1(\bar{t}_1)$ stopped on hitting $\overline{\mathbb{H} \setminus H_2}$.

## 6. Constructing new martingales.

THEOREM 6.1. *For any $(H_1^m, H_2^m) \in \text{HP}$, $1 \leq m \leq n$, there is a continuous function $M_*(t_1, t_2)$ defined on $[0, \infty)^2$ that satisfies the following properties:* (i) $M_* = M$ *on* $[0, T_1(H_1^m)] \times [0, T_2(H_2^m)]$ *for* $m = 1, \ldots, n$; (ii) $M_*(t, 0) = M_*(0, t) = 1$ *for any $t \geq 0$;* (iii) $M_*(t_1, t_2) \in [C_1, C_2]$ *for any $t_1, t_2 \geq 0$, where $C_2 > C_1 > 0$ are constants depending only on $H_j^m$, $j = 1, 2$, $1 \leq m \leq n$;* (iv) *for any $(\mathcal{F}_t^2)$-stopping time $\bar{t}_2$, $(M_*(t_1, \bar{t}_2), t_1 \geq 0)$ is a bounded continuous $(\mathcal{F}_{t_1}^1 \times \mathcal{F}_{\bar{t}_2}^2)_{t_1 \geq 0}$-martingale; and* (v) *for any $(\mathcal{F}_t^1)$-stopping time $\bar{t}_1$, $(M_*(\bar{t}_1, t_2), t_2 \geq 0)$ is a bounded continuous $(\mathcal{F}_{\bar{t}_1}^1 \times \mathcal{F}_{t_2}^2)_{t_2 \geq 0}$-martingale.*

PROOF. We will first define $M_*$ and then check its properties. The first quadrant $[0, \infty)^2$ is divided by the horizontal or vertical lines $\{x_j = T_j(H_j^m)\}$, $1 \leq m \leq n$, $j = 1, 2$, into small rectangles, and $M_*$ is piecewise defined on each rectangle. Theorem 4.2 will be used to prove the martingale properties.

Let $\mathbb{N}_n := \{k \in \mathbb{N} : k \leq n\}$. Write $T_j^k$ for $T_j(H_j^k)$, $k \in \mathbb{N}_n$, $j = 1, 2$. Let $S \subset \mathbb{N}_n$ be such that $\bigcup_{k \in S} [0, T_1^k] \times [0, T_2^k] = \bigcup_{k=1}^n [0, T_1^k] \times [0, T_2^k]$, and $\sum_{k \in S} k \leq \sum_{k \in S'} k$ if $S' \subset \mathbb{N}_n$ also satisfies this property. Such $S$ is a random nonempty set, and $|S| \in \mathbb{N}_n$ is a random number. Define a partial order "$\preceq$" on $[0, \infty]^2$ such that $(s_1, s_2) \preceq (t_1, t_2)$ iff $s_1 \leq t_1$ and $s_2 \leq t_2$. If $(s_1, s_2) \preceq (t_1, t_2)$ and $(s_t, s_2) \neq (t_1, t_2)$, we write $(s_1, s_2) \prec (t_1, t_2)$. Then for each $m \in \mathbb{N}_n$, there is $k \in S$ such that $(T_1^m, T_2^m) \preceq (T_1^k, T_2^k)$; and for each $k \in S$, there is no $m \in \mathbb{N}_n$ such that $(T_1^k, T_2^k) \prec (T_1^m, T_2^m)$.

There is a map $\sigma$ from $\{1, \ldots, |S|\}$ onto $S$ such that if $1 \leq k_1 < k_2 \leq |S|$, then

$$(6.1) \qquad T_1^{\sigma(k_1)} < T_1^{\sigma(k_2)}, \qquad T_2^{\sigma(k_1)} > T_2^{\sigma(k_2)}.$$

Define $T_1^{\sigma(0)} = T_2^{\sigma(|S|+1)} = 0$ and $T_1^{\sigma(|S|+1)} = T_2^{\sigma(0)} = \infty$. Then (6.1) still holds for $0 \le k_1 < k_2 \le |S|+1$.

Extend the definition of $M$ to $[0,\infty] \times \{0\} \cup \{0\} \times [0,\infty]$ such that $M(t,0) = M(0,t) = 1$ for $t \ge 0$. Fix $(t_1, t_2) \in [0,\infty]^2$. There are $k_1 \in \mathbb{N}_{|S|+1}$ and $k_2 \in \mathbb{N}_{|S|} \cup \{0\}$ such that

(6.2) $$T_1^{\sigma(k_1-1)} \le t_1 \le T_1^{\sigma(k_1)}, \qquad T_2^{\sigma(k_2+1)} \le t_2 \le T_2^{\sigma(k_2)}.$$

If $k_1 \le k_2$, let

(6.3) $$M_*(t_1, t_2) = M(t_1, t_2).$$

It $k_1 \ge k_2 + 1$, let

(6.4) $$\begin{aligned}
M_*(t_1, t_2) = &(M(T_1^{\sigma(k_2)}, t_2) M(T_1^{\sigma(k_2+1)}, T_2^{\sigma(k_2+1)}) \\
&\cdots M(T_1^{\sigma(k_1-1)}, T_2^{\sigma(k_1-1)}) M(t_1, T_2^{\sigma(k_1)})) \\
&\times (M(T_1^{\sigma(k_2)}, T_2^{\sigma(k_2+1)}) \\
&\cdots M(T_1^{\sigma(k_1-2)}, T_2^{\sigma(k_1-1)}) M(T_1^{\sigma(k_1-1)}, T_2^{\sigma(k_1)}))^{-1}.
\end{aligned}$$

In the above formula, there are $k_1 - k_2 + 1$ terms in the numerator, and $k_1 - k_2$ terms in the denominator. For example, if $k_1 - k_2 = 1$, then

$$M_*(t_1, t_2) = M(T_1^{\sigma(k_2)}, t_2) M(t_1, T_2^{\sigma(k_1)}) / M(T_1^{\sigma(k_2)}, T_2^{\sigma(k_1)}).$$

We need to show that $M_*(t_1, t_2)$ is well defined. First, we show that the $M(\cdot, \cdot)$ in (6.3) and (6.4) are defined. Note that $M$ is defined on

$$Z := \bigcup_{k=0}^{|S|+1} [0, T_1^{\sigma(k)}] \times [0, T_2^{\sigma(k)}].$$

If $k_1 \le k_2$, then $t_1 \le T_1^{\sigma(k_1)} \le T_1^{\sigma(k_2)}$ and $t_2 \le T_2^{\sigma(k_2)}$, so $(t_1, t_2) \in Z$. Thus, $M(t_1, t_2)$ in (6.3) is defined. Now suppose $k_1 \ge k_2 + 1$. Since $t_2 \le T_2^{\sigma(k_2)}$ and $t_1 \le T_1^{\sigma(k_1)}$, so $(T_1^{\sigma(k_2)}, t_2), (t_1, T_2^{\sigma(k_1)}) \in Z$. It is clear that $(T_1^{\sigma(k)}, T_2^{\sigma(k)}) \in Z$ for $k_2 + 1 \le k \le k_1 - 1$. Thus, the $M(\cdot, \cdot)$ in the numerator of (6.4) are defined. For $k_2 \le k \le k_1 - 1$, $T_1^{\sigma(k)} \le T_1^{\sigma(k+1)}$, so $(T_1^{\sigma(k)}, T_2^{\sigma(k+1)}) \in Z$. Thus, the $M(\cdot, \cdot)$ in the denominator of (6.4) are defined.

Second, we show that the value of $M_*(t_1, t_2)$ does not depend on the choice of $(k_1, k_2)$ that satisfies (6.2). Suppose (6.2) holds with $(k_1, k_2)$ replaced by $(k_1', k_2)$, and $k_1' \ne k_1$. Then $|k_1' - k_1| = 1$. We may assume $k_1' = k_1 + 1$. Then $t_1 = T_1^{\sigma(k_1)}$. Let $M_*'(t_1, t_2)$ denote the $M_*(t_1, t_2)$ defined using $(k_1', k_2)$. There are three cases.

*Case* 1. $k_1 < k_1' \le k_2$. Then from (6.3), $M_*'(t_1, t_2) = M(t_1, t_2) = M_*(t_1, t_2)$.

*Case* 2. $k_1 = k_2$ and $k'_1 - k_2 = 1$. Then $T_1^{\sigma(k_2)} = T_1^{\sigma(k_1)} = t_1$. So from (6.3) and (6.4),

$$M'_*(t_1, t_2) = M(T_1^{\sigma(k_2)}, t_2) M(t_1, T_2^{\sigma(k_1)}) / M(T_1^{\sigma(k_2)}, T_2^{\sigma(k_1)})$$
$$= M(t_1, t_2) = M_*(t_1, t_2).$$

*Case* 3. $k'_1 > k_1 > k_2$. From (6.4) and that $T_1^{\sigma(k_1)} = t_1$, we have

$M'_*(t_1, t_2)$

$$= \frac{M(T_1^{\sigma(k_2)}, t_2) M(T_1^{\sigma(k_2+1)}, T_2^{\sigma(k_2+1)}) \cdots M(T_1^{\sigma(k_1)}, T_2^{\sigma(k_1)}) M(t_1, T_2^{\sigma(k_1+1)})}{M(T_1^{\sigma(k_2)}, T_2^{\sigma(k_2+1)}) \cdots M(T_1^{\sigma(k_1-1)}, T_2^{\sigma(k_1)}) M(T_1^{\sigma(k_1)}, T_2^{\sigma(k_1+1)})}$$

$$= \frac{M(T_1^{\sigma(k_2)}, t_2) M(T_1^{\sigma(k_2+1)}, T_2^{\sigma(k_2+1)}) \cdots M(t_1, T_2^{\sigma(k_1)})}{M(T_1^{\sigma(k_2)}, T_2^{\sigma(k_2+1)}) \cdots M(T_1^{\sigma(k_1-1)}, T_2^{\sigma(k_1)})} = M_*(t_1, t_2).$$

Similarly, if (6.2) holds with $(k_1, k_2)$ replaced by $(k_1, k'_2)$, then $M_*(t_1, t_2)$ defined using $(k_1, k'_2)$ has the same value as $M(t_1, t_2)$. Thus, $M_*$ is well defined.

From the definition, it is clear that for each $k_1 \in \mathbb{N}_{|S|+1}$ and $k_2 \in \mathbb{N}_{|S|} \cup \{0\}$, $M_*$ is continuous on $[T_1^{\sigma(k_1-1)}, T_1^{\sigma(k_1)}] \times [T_2^{\sigma(k_2+1)}, T_1^{\sigma(k_2)}]$. Thus, $M_*$ is continuous on $[0, \infty)^2$. Let $(t_1, t_2) \in [0, \infty)^2$. Suppose $(t_1, t_2) \in [0, T_1^m] \times [0, T_2^m]$ for some $m \in \mathbb{N}_n$. There is $k \in \mathbb{N}_{|S|}$ such that $(T_1^m, T_2^m) \preceq (T_1^{\sigma(k)}, T_2^{\sigma(k)})$. Then we may choose $k_1 \leq k$ and $k_2 \geq k$ such that (6.2) holds, so $M_*(t_1, t_2) = M(t_1, t_2)$. Thus, (i) is satisfied. If $t_1 = 0$, we may choose $k_1 = 1$ in (6.2). Then either $k_1 \leq k_2$ or $k_2 = 0$. If $k_1 \leq k_2$, then $M_*(t_1, t_2) = M(t_1, t_2) = 1$ because $t_1 = 0$. If $k_2 = 0$, then

$$M_*(t_1, t_2) = M(T_1^{\sigma(0)}, t_2) M(t_1, T_2^{\sigma(1)}) / M(T_1^{\sigma(0)}, T_2^{\sigma(1)}) = 1$$

because $T_1^{\sigma(0)} = t_1 = 0$. Similarly, $M_*(t_1, t_2) = 0$ if $t_2 = 0$. So (ii) is also satisfied. And (iii) follows from Lemma 5.1 and the definition of $M_*$.

Now we prove (iv). Suppose $(t_1, t_2) \in [0, \infty)^2$ and $t_2 \geq \bigvee_{m=1}^n T_2^m = T_2^{\sigma(1)}$. Then (6.2) holds with $k_2 = 0$ and some $k_1 \in \{1, \ldots, |S|+1\}$. So $k_1 \geq k_2 + 1$. Since $T_1^{\sigma(k_2)} = 0$ and $M(0, t) = 1$ for any $t \geq 0$, so from (6.4) we have

$$M_*(t_1, t_2) = \frac{M(T_1^{\sigma(k_2+1)}, T_2^{\sigma(k_2+1)}) \cdots M(T_1^{\sigma(k_1-1)}, T_2^{\sigma(k_1-1)}) M(t_1, T_2^{\sigma(k_1)})}{M(T_1^{\sigma(k_2+1)}, T_2^{\sigma(k_2+2)}) \cdots M(T_1^{\sigma(k_1-1)}, T_2^{\sigma(k_1)})}.$$

The right-hand side of the above equality has no $t_2$. So $M_*(t_1, t_2) = M_*(t_1, \bigvee_{m=1}^n T_2^m)$ for any $t_2 \geq \bigvee_{m=1}^n T_2^m$. Similarly, $M_*(t_1, t_2) = M_*(\bigvee_{m=1}^n T_1^m, t_2)$ for any $t_1 \geq \bigvee_{m=1}^n T_1^m$.

Fix an $(\mathcal{F}_t^2)$-stopping time $\bar{t}_2$. Since $M_*(\cdot, \bar{t}_2) = M_*(\cdot, \bar{t}_2 \wedge (\bigvee_{m=1}^n T_2^m))$, and $\bar{t}_2 \wedge (\bigvee_{m=1}^n T_2^m)$ is also an $(\mathcal{F}_t^2)$-stopping time, so we may assume that $\bar{t}_2 \leq \bigvee_{m=1}^n T_2^m$. Let $I_0 = \bar{t}_2$. For $s \in \mathbb{N} \cup \{0\}$, define

$$R_s = \sup\{T_1^m : m \in \mathbb{N}_n, T_2^m \geq I_s\};$$
(6.5)
$$I_{s+1} = \sup\{T_2^m : m \in \mathbb{N}_n, T_2^m < I_s, T_1^m > R_s\}.$$

Here we set $\sup(\varnothing) = 0$. Then we have a nondecreasing sequence $(R_s)$ and a non-increasing sequence $(I_s)$. Let $S$ and $\sigma(k)$, $0 \leq k \leq |S|+1$, be as in the definition of $M_*$. From the property of $S$, for any $s \in \mathbb{N} \cup \{0\}$,

(6.6) $$R_s = \sup\{T_1^k : k \in S, T_2^k \geq I_s\}.$$

Suppose for some $s \in \mathbb{N} \cup \{0\}$, there is $m \in \mathbb{N}_n$ that satisfies $T_2^m < I_s$ and $T_1^m > R_s$. Then there is $k \in S$ such that $T_j^k \geq T_j^m$, $j = 1, 2$. If $T_2^k \geq I_s$, then from (6.6) we have $R_s \geq T_1^k \geq T_1^m$, which contradicts that $T_1^m > R_s$. Thus, $T_2^k < I_s$. Now $T_2^k < I_s$, $T_1^k \geq T_1^m > R_s$, and $T_2^k \geq T_2^m$. Thus, for any $s \in \mathbb{N} \cup \{0\}$,

(6.7) $$I_{s+1} = \sup\{T_2^k : k \in S, T_2^k < I_s, T_1^k > R_s\}.$$

First suppose $\bar{t}_2 > 0$. Since $\bar{t}_2 \leq \bigvee_{m=1}^n T_2^m = T_2^{\sigma(0)}$, so there is a unique $k_2 \in \mathbb{N}_{|S|}$ such that $T_2^{\sigma(k_2)} \geq \bar{t}_2 > T_2^{\sigma(k_2+1)}$. From (6.6) and (6.7), we have $R_s = T_1^{\sigma(k_2+s)}$ for $0 \leq s \leq |S| - k_2$; $R_s = T_1^{\sigma(|S|)}$ for $s \geq |S| - k_2$; $I_s = T_2^{\sigma(k_2+s)}$ for $1 \leq s \leq |S| - k_2$; and $I_s = 0$ for $s \geq |S| - k_2 + 1$. Since $R_0 = T_1^{\sigma(k_2)}$ and $\bar{t}_2 \leq T_2^{\sigma(k_2)}$, so from (i),

(6.8) $$M_*(t_1, \bar{t}_2) = M(t_1, \bar{t}_2) \qquad \text{for } t_1 \in [0, R_0].$$

Suppose $t_1 \in [R_{s-1}, R_s]$ for some $s \in \mathbb{N}_{|S|-k_2}$. Let $k_1 = k_2 + s$. Then $T_1^{\sigma(k_1-1)} \leq t_1 \leq T_1^{\sigma(k_1)}$. Since $I_s = T_2^{\sigma(k_2+s)} = T_2^{\sigma(k_1)}$, so from (6.4),

$$M_*(t_1, \bar{t}_2)/M_*(R_{s-1}, \bar{t}_2) = M(t_1, I_s)/M(R_{s-1}, I_s),$$
(6.9)
$$\text{for } t_1 \in [R_{s-1}, R_s].$$

Note that if $s \geq |S| - k_2 + 1$, (6.9) still holds because $R_s = R_{s-1}$. Suppose $t_1 \geq R_n$. Since $n \geq |S| - k_2$, so $R_n = T_1^{\sigma(|S|)} = \bigvee_{m=1}^n T_1^m$. From the discussion at the beginning of the proof of (iv), we have

(6.10) $$M_*(t_1, \bar{t}_2) = M_*(R_n, \bar{t}_2), \qquad \text{for } t_1 \in [R_n, \infty].$$

If $\bar{t}_2 = 0$, (6.8)–(6.10) still hold because all $I_s = 0$ and so $M_*(t_1, \bar{t}_2) = M(t_1, I_s) = M(t_1, 0) = 1$ for any $t_1 \geq 0$.

Let $R_{-1} = 0$. We claim that for each $s \in \mathbb{N} \cup \{0\}$, $R_s$ is an $(\mathcal{F}_t^1 \times \mathcal{F}_{\bar{t}_2}^2)_{t \geq 0}$-stopping time and $I_s$ is $\mathcal{F}_{R_{s-1}, \bar{t}_2}$-measurable. Recall that $\mathcal{F}_{R_{s-1}, \bar{t}_2}$ is the $\sigma$-algebra obtained from the filtration $(\mathcal{F}_t^1 \times \mathcal{F}_{\bar{t}_2}^2)_{t \geq 0}$ and its stopping time $R_{s-1}$. It is clear that $R_{-1} = 0$ is an $(\mathcal{F}_t^1 \times \mathcal{F}_{\bar{t}_2}^2)_{t\geq 0}$-stopping time, and $I_0 = \bar{t}_2$ is $\mathcal{F}_{R_{-1},\bar{t}_2}$-measurable. Now suppose $I_s$ is $\mathcal{F}_{R_{s-1},\bar{t}_2}$-measurable. Since $I_s \leq \bar{t}_2$ and $R_{s-1} \leq R_s$, so for any $t \geq 0$, $\{R_s \leq t\} = \{R_{s-1} \leq t\} \cap \mathcal{E}_t$, where

$$\mathcal{E}_t = \bigcap_{m=1}^n (\{T_2^m < I_s\} \cup \{T_1^m \leq t\})$$

$$= \bigcap_{m=1}^n \left( \bigcup_{q \in \mathbb{Q}} (\{T_2^m < q \leq \bar{t}_2\} \cap \{q < I_s\}) \cup \{T_1^m \leq t\} \right).$$

Thus, $\mathcal{E}_t \in \mathcal{F}_{R_{s-1},\bar{t}_2} \vee (\mathcal{F}_t^1 \times \mathcal{F}_{\bar{t}_2}^2)$, and so $\{R_s \leq t\} \in \mathcal{F}_t^1 \times \mathcal{F}_{\bar{t}_2}^2$ for any $t \geq 0$. Therefore, $R_s$ is an $(\mathcal{F}_t^1 \times \mathcal{F}_{\bar{t}_2}^2)_{t\geq 0}$-stopping time. Next we consider $I_{s+1}$. For any $h \geq 0$,

$$\{I_{s+1} > h\} = \bigcup_{m=1}^n (\{h < T_2^m < I_s\} \cap \{T_1^m > R_s\})$$

$$= \bigcup_{m=1}^n \left( \bigcup_{q \in \mathbb{Q}} (\{h < T_2^m < q < \bar{t}_2\} \cap \{q < I_s\}) \cap \{T_1^m > R_s\} \right) \in \mathcal{F}_{R_s,\bar{t}_2}.$$

Thus, $I_{s+1}$ is $\mathcal{F}_{R_s,\bar{t}_2}$-measurable. So the claim is proved by induction.

Since $\bar{t}_2 \leq \bigvee_{m=1}^n T_2^m < T_2$, so from Theorem 4.2, for any $s \in \mathbb{N}_n$, $(M(R_{s-1} + t, I_s), 0 \leq t < T_1(I_s) - R_{s-1})$ is a continuous $(\mathcal{F}_{R_{s-1}+t,\bar{t}_2})_{t\geq 0}$-local martingale. For $m \in \mathbb{N}_n$, if $T_2^m \geq I_s$, then $T_1^m < T_1(T_2^m) \leq T_1(I_s)$. So from (6.5) we have $R_s < T_1(I_s)$. From (6.9), we find that $(M_*(R_{s-1} + t, \bar{t}_2), 0 \leq t \leq R_s - R_{s-1})$ is a continuous $(\mathcal{F}_{R_{s-1}+t,\bar{t}_2})_{t\geq 0}$-local martingale for any $s \in \mathbb{N}_n$. From Theorem 4.1 and (6.8), $(M_*(t, \bar{t}_2), 0 \leq t \leq R_0)$ is a continuous $(\mathcal{F}_{t,\bar{t}_2})_{t\geq 0}$-local martingale. From (6.10), $(M_*(R_n + t, \bar{t}_2), t \geq 0)$ is a continuous $(\mathcal{F}_{R_n+t,\bar{t}_2})_{t\geq 0}$-local martingale. Thus, $(M_*(t, \bar{t}_2), t \geq 0)$ is a continuous $(\mathcal{F}_{t,\bar{t}_2})_{t\geq 0}$-local martingale. Since by (iii) $M_*(t_1, t_2) \in [C_1, C_2]$, so this local martingale is a bounded martingale. Thus, (iv) is satisfied. Finally, (v) follows from the symmetry in the definition of (6.3) and (6.4) of $M_*$. □

## 7. Coupling measures.

PROOF OF THEOREM 2.1. From conformal invariance, we may assume that $D = \mathbb{H}$, $a = x_1$ and $b = x_2$. Let $\xi_j(t)$ and $\beta_j(t)$, $0 \leq t \leq T_j$, $j = 1, 2$, be as in Section 4. For $j = 1, 2$, let $\mu_j$ denote the distribution of $(\xi_j(t), 0 \leq t \leq T_j)$. Let $\mu = \mu_1 \times \mu_2$. Then $\mu$ is the joint distribution of $\xi_1$ and $\xi_2$, since they are independent.

Let $\widehat{\mathbb{C}} = \mathbb{C} \cup \{\infty\}$ be the Riemann sphere with spherical metric. Let $\Gamma_{\widehat{\mathbb{C}}}$ denote the space of nonempty compact subsets of $\widehat{\mathbb{C}}$ endowed with the Hausdorff metric. Then $\Gamma_{\widehat{\mathbb{C}}}$ is a compact metric space. For a chordal Loewner trace $\beta(t)$, $0 \leq t \leq T$, let $G(\beta) := \{\beta(t) : 0 \leq t \leq T\} \in \Gamma_{\widehat{\mathbb{C}}}$. For $j = 1, 2$, let $\bar{\mu}_j$ denote the distribution of $G(\beta_j)$, which is a probability measure on $\Gamma_{\widehat{\mathbb{C}}}$. We want to prove that $\bar{\mu}_1 = \bar{\mu}_2$. Let $\bar{\mu} = \bar{\mu}_1 \times \bar{\mu}_2$, which is the joint distribution of $G(\beta_1)$ and $G(\beta_2)$.

Let $\mathrm{HP}_*$ be the set of $(H_1, H_2) \in \mathrm{HP}$ such that, for $j = 1, 2$, $H_j$ is a polygon whose vertices have rational coordinates. Then $\mathrm{HP}_*$ is countable. Let $(H_1^m, H_2^m)$, $m \in \mathbb{N}$, be an enumeration of $\mathrm{HP}_*$. For each $n \in \mathbb{N}$, let $M_*^n(t_1, t_2)$ be the $M_*(t_1, t_2)$ given by Theorem 6.1 for $(H_1^m, H_2^m)$, $1 \leq m \leq n$, in the above enumeration.

For each $n \in \mathbb{N}$, define $\nu^n = (\nu_1^n, \nu_2^n)$ such that $d\nu^n/d\mu = M_*^n(\infty, \infty)$. From Theorem 6.1, $M_*^n(\infty, \infty) > 0$ and $\int M_*^n(\infty, \infty) \, d\mu = \mathbf{E}[M_*^n(\infty, \infty)] = 1$, so $\nu^n$ is a probability measure. Then $d\nu_1^n/d\mu_1 = \mathbf{E}[M_*^n(\infty, \infty)|\mathcal{F}_\infty^2] = M_*^n(\infty, 0) = 1$. Thus, $\nu_1^n = \mu_1$. Similarly, $\nu_2^n = \mu_2$. So each $\nu^n$ is a coupling of $\mu_1$ and $\mu_2$.

For each $n \in \mathbb{N}$, suppose $(\zeta_1^n(t), 0 \leq t \leq S_1^n)$ and $(\zeta_2^n(t), 0 \leq t \leq S_2^n)$ have the joint distribution $\nu^n$. Let $\gamma_j^n(t)$, $0 \leq t \leq S_j$, $j = 1, 2$, be the chordal Loewner trace driven by $\zeta_j^n$. Let $\bar{\nu}^n = (\bar{\nu}_1^n, \bar{\nu}_2^n)$ denote the joint distribution of $G(\gamma_1^n)$ and $G(\gamma_2^n)$. Since $\Gamma_{\widehat{\mathbb{C}}} \times \Gamma_{\widehat{\mathbb{C}}}$ is compact, so $(\bar{\nu}^n, n \in \mathbb{N})$ has a subsequence $(\bar{\nu}^{n_k} : k \in \mathbb{N})$ that converges weakly to some probability measure $\bar{\nu} = (\bar{\nu}_1, \bar{\nu}_2)$ on $\Gamma_{\widehat{\mathbb{C}}} \times \Gamma_{\widehat{\mathbb{C}}}$. Then for $j = 1, 2$, $\bar{\nu}_j^{n_k} \to \bar{\nu}_j$ weakly. For $n \in \mathbb{N}$ and $j = 1, 2$, since $\nu_j^n = \mu_j$, so $\bar{\nu}_j^n = \bar{\mu}_j$. Thus, $\bar{\nu}_j = \bar{\mu}_j$, $j = 1, 2$. So $\bar{\nu}_j$, $j = 1, 2$, is supported by the space of graphs of crosscuts in $\mathbb{H}$. From Proposition 2.2, there are $\zeta_1 \in C([0, S_1])$ and $\zeta_2 \in C([0, S_2])$ such that the joint distribution of $G(\gamma_1)$ and $G(\gamma_2)$ is $\bar{\nu}$, where $\gamma_j(t)$ is the chordal Loewner trace driven by $\zeta_j(t)$, $j = 1, 2$.

Now fix $m \in \mathbb{N}$. From Theorem 4.1, $M(T_1(H_1^m), T_2(H_2^m))$ is positive and $\mathcal{F}^1_{T_1(H_1^m)} \times \mathcal{F}^2_{T_2(H_2^m)}$-measurable, and $\int M(T_1(H_1^m), T_2(H_2^m)) \, d\mu = 1$. Define $\nu_{(m)}$ on $\mathcal{F}^1_{T_1(H_1^m)} \times \mathcal{F}^2_{T_2(H_2^m)}$ such that $d\nu_{(m)}/d\mu = M(T_1(H_1^m), T_2(H_2^m))$. Then $\nu_{(m)}$ is a probability measure. From Theorem 6.1, if $n \geq m$, then

$$\frac{d\nu^n}{d\mu}\bigg|_{\mathcal{F}^1_{T_1(H_1^m)} \times \mathcal{F}^2_{T_2(H_2^m)}} = \mathbf{E}[M_*^n(\infty, \infty)|\mathcal{F}^1_{T_1(H_1^m)} \times \mathcal{F}^2_{T_2(H_2^m)}]$$

$$= M_*^n(T_1(H_1^m), T_2(H_2^m)) = M(T_1(H_1^m), T_2(H_2^m)).$$

Thus, $\nu_{(m)}$ equals the restriction of $\nu^n$ to $\mathcal{F}^1_{T_1(H_1^m)} \times \mathcal{F}^2_{T_2(H_2^m)}$ if $n \geq m$.

For a chordal Loewner trace $\gamma(t)$, $0 \leq t \leq S$, and a hull $H$ in $\mathbb{H}$ w.r.t. 0 that contains some neighborhood of $\gamma(0)$ in $\mathbb{H}$, let $G_H(\gamma) := \{\gamma(t) : 0 \leq t \leq T_H\} \in \Gamma_{\widehat{\mathbb{C}}}$, where $T_H$ is the first $t$ such that $\gamma(t) \in \overline{\mathbb{H} \setminus H}$ or $t = S$. Then $G_H(\gamma) \subset G(\gamma)$. Let $\bar{\nu}^n_{(m)}$ denote the distribution of $(G_{H_1^m}(\gamma_1^n), G_{H_2^m}(\gamma_2^n))$. Then $\bar{\nu}^n_{(m)}$ is determined by the distribution of $(\zeta_1^n, \zeta_2^n)$ restricted to $\mathcal{F}^1_{T_1(H_1^m)} \times \mathcal{F}^2_{T_2(H_2^m)}$, which equals $\nu_{(m)}$ if $n \geq m$. Let $\bar{\nu}_{(m)} = \bar{\nu}^m_{(m)}$. Then $\bar{\nu}^n_{(m)} = \bar{\nu}_{(m)}$ for $n \geq m$.

Let $\tau_{(m)}^{n_k}$ denote the distribution of $(G(\gamma_1^{n_k}), G(\gamma_2^{n_k}), G_{H_1^m}(\gamma_1^{n_k}), G_{H_2^m}(\gamma_2^{n_k}))$. Then $\tau_{(m)}^{n_k}$ is supported by $\Xi$, which is the set of $(L_1, L_2, F_1, F_2) \in \Gamma_{\widehat{\mathbb{C}}}^4$ such that $F_j \subset L_j$ for $j = 1, 2$. It is easy to check that $\Xi$ is a closed subset of $\Gamma_{\widehat{\mathbb{C}}}^4$. Then $(n_k)$ has a subsequence $(n_k')$ such that $\tau_{(m)}^{n_k'}$ converges weakly to some probability measure $\tau_{(m)}$ on $\Xi$. Since the marginal of $\tau_{(m)}^{n_k'}$ at the first two variables equals $\bar\nu^{n_k'}$, and $\bar\nu^{n_k'} \to \bar\nu$ weakly, so the marginal of $\tau_{(m)}$ at the first two variables equals $\bar\nu$. Since the marginal of $\tau_{(m)}^{n_k'}$ at the last two variables equals $\bar\nu_{(m)}^{n_k'}$, which equals $\bar\nu_{(m)}$ if $n_k' \geq m$, so the marginal of $\tau_{(m)}$ at the last two variables equals $\bar\nu_{(m)}$.

Let the $\Xi$-valued random variable $(L_1, L_2, F_1, F_2)$ have the distribution $\tau_{(m)}$. Then $\bar\nu$ is the distribution of $(L_1, L_2)$ and $\bar\nu_{(m)}$ is the distribution of $(F_1, F_2)$. Note that $\bar\nu_{(m)}$ is supported by the space of pairs of curves $(\alpha_1, \alpha_2)$ such that, for $j = 1, 2$, $\alpha_j$ is a simple curve whose one end is $x_j$, the other end lies on $\partial H_j^m \cap \mathbb{H}$, and whose other part lies in the interior of $H_j^m$. For $j = 1, 2$, since $L_j = G(\gamma_j)$, so from the properties of $\Xi$ and $\bar\nu_{(m)}$, we have $F_j = G_{H_j^m}(\gamma_j)$, which means that $(G_{H_1^m}(\gamma_1), G_{H_2^m}(\gamma_2))$ has the distribution $\bar\nu_{(m)}$. Since the distribution of $(G_{H_1^m}(\gamma_1), G_{H_2^m}(\gamma_2))$ determines the distribution of $(\zeta_1, \zeta_2)$ restricted to $\mathcal{F}_{T_1(H_1^m)}^1 \times \mathcal{F}_{T_2(H_2^m)}^2$, so the distribution of $(\zeta_1, \zeta_2)$ restricted to $\mathcal{F}_{T_1(H_1^m)}^1 \times \mathcal{F}_{T_2(H_2^m)}^2$ equals $\nu_{(m)}$. Since $d\nu_{(m)}/d\mu = M(T_1(H_1^m), T_2(H_2^m))$, so from the discussion after the proof of Theorem 5.1, for any $(\mathcal{F}_t^2)$-stopping time $\bar t_2$ with $\bar t_2 \leq T_2(H_2^m)$, $(\gamma_1(t), 0 \leq t \leq T_1(H_1^m))$ is a time-change of a chordal SLE$_\kappa$ trace in $\mathbb{H} \setminus \gamma_2((0, \bar t_2])$ from $x_1$ to $\gamma_2(\bar t_2)$ stopped on hitting $\overline{\mathbb{H} \setminus H_1^m}$.

Now fix an $(\mathcal{F}_t^2)$-stopping time $\bar t_2$ with $\bar t_2 < T_2$. Recall that $T_1(\bar t_2)$ is the maximal such that $\gamma_1([0, T_1(\bar t_2)))$ is disjoint from $\gamma_2([0, \bar t_2])$. For $n \in \mathbb{N}$, define

$$R_n = \sup\{T_1(H_1^m) : m \in \mathbb{N}_n, \bar t_2 \leq T_2(H_2^m)\}.$$

Here we set $\sup(\varnothing) = 0$. Then for any $t \geq 0$,

$$\{R_n \leq t\} = \bigcap_{m=1}^n (\{\bar t_2 > T_2(H_2^m)\} \cup \{T_1(H_1^m) \leq t\}) \in \mathcal{F}_t^1 \times \mathcal{F}_{\bar t_2}^2.$$

So $R_n$ is an $(\mathcal{F}_t^1 \times \mathcal{F}_{\bar t_2}^2)_{t \geq 0}$-stopping time for each $n \in \mathbb{N}$. For $m \in \mathbb{N}_n$, let $\bar t_2^m = \bar t_2 \wedge T_2(H_2^m)$. Then $\bar t_2^m$ is an $(\mathcal{F}_t^2)$-stopping time, and $\bar t_2^m \leq T_2(H_2^m)$. From the last paragraph, we conclude that $\gamma_1(t)$, $0 \leq t \leq T_1(H_1^m)$, is a time-change of a part of the chordal SLE$_\kappa$ trace in $\mathbb{H} \setminus \gamma_1((0, \bar t_2^m])$ from $x_1$ to $\gamma_2(\bar t_2^m)$. Let $\mathcal{E}_{n,m} = \{\bar t_2 \leq T_2(H_2^m)\} \cap \{R_n = T_1(H_1^m)\}$. Since on each $\mathcal{E}_{n,m}$, $\bar t_2 = \bar t_2^m$ and $R_n = T_1(H_1^m)$, and $\{R_n > 0\} = \bigcup_{m=1}^n \mathcal{E}_{n,m}$, so $\gamma_1(t)$, $0 \leq t \leq R_n$, is a time-change of a part of the chordal SLE$_\kappa$ trace in $\mathbb{H} \setminus \gamma_1((0, \bar t_2])$ from $x_1$ to $\gamma_2(\bar{t}_2)$. Let $R_\infty = \bigvee_{n=1}^\infty R_n$. Then $\gamma_1(t)$, $0 \leq t < R_\infty$, is a time-change of a part of the chordal $\text{SLE}_\kappa$ trace in $\mathbb{H} \setminus \gamma_1((0, \bar{t}_2])$ from $x_1$ to $\gamma_2(\bar{t}_2)$.

For each $n \in \mathbb{N}$ and $m \in \mathbb{N}_n$, if $\bar{t}_2 \leq T_2(H_2^m)$, then $T_1(H_2^m) < T_1(\bar{t}_2)$, so $R_n < T_1(\bar{t}_2)$. Thus, $R_\infty \leq T_1(\bar{t}_2)$. If $R_\infty < T_1(\bar{t}_2)$, then $\gamma_1((0, R_\infty])$ is disjoint from $\gamma_2((0, \bar{t}_2])$, so there is $(H_1^m, H_2^m) \in \text{HP}_*$ such that $\gamma_1((0, R_\infty])$ and $\gamma_2((0, \bar{t}_2])$ are contained in the interiors of $H_1^m$ and $H_2^m$, respectively. Then $\bar{t}_2 \leq T_2(H_2^m)$ and $R_m \leq R_\infty < T_1(H_1^m)$, which contradicts the definition of $R_m$. Thus, $R_\infty = T_1(\bar{t}_2)$. So $\gamma_1(t)$, $0 \leq t < T_1(\bar{t}_2)$, is a time-change of a part of the chordal $\text{SLE}_\kappa$ trace in $\mathbb{H} \setminus \gamma_1((0, \bar{t}_2])$ from $x_1$ to $\gamma_2(\bar{t}_2)$. From the definition of $T_1(\bar{t}_2)$ we have $\gamma_1(T_1(\bar{t}_2)) \in G(\gamma_2)$. Thus, $\gamma_1(t)$, $0 \leq t < T_1(\bar{t}_2)$, is a time-change of a full chordal $\text{SLE}_\kappa$ trace in $\mathbb{H} \setminus \gamma_1((0, \bar{t}_2])$ from $x_1$ to $\gamma_2(\bar{t}_2)$. Since $\kappa \in (0, 4]$, so almost surely $\gamma_1(T_1(\bar{t}_2)) = \gamma_2(\bar{t}_2)$. Thus, $\gamma_2(\bar{t}_2) \in G(\gamma_1)$ almost surely.

For $n \in \mathbb{N}$ and $q \in \mathbb{Q}_{\geq 0}$, let $\bar{t}_2^{n,q} = q \wedge T_2(H_2^n)$. Then each $\bar{t}_2^{n,q}$ is an $(\mathcal{F}_t^2)$-stopping time with $\bar{t}_2^{q,n} < T_2$. Since $\mathbb{N} \times \mathbb{Q}_{\geq 0}$ is countable, so almost surely $\gamma_2(\bar{t}_2^{q,n}) \in G(\gamma_1)$ for every $n \in \mathbb{N}$ and $q \in \mathbb{Q}_{\geq 0}$. Since $\mathbb{Q}_{\geq 0}$ is dense in $\mathbb{R}_{\geq 0}$, $\gamma_2$ is continuous, and $G(\gamma_1)$ is closed, so almost surely for every $n \in \mathbb{N}$, $\gamma_2([0, T_2(H_2^n)]) \subset G(\gamma_1)$. Since $T_2 = \bigvee_{n=1}^\infty T_2(H_2^n)$, so $G(\gamma_2) \subset G(\gamma_1)$ almost surely. Similarly, $G(\gamma_1) \subset G(\gamma_2)$ almost surely. Thus, $G(\gamma_1) = G(\gamma_2)$ almost surely. Since for $j = 1, 2$, the distribution of $G(\gamma_j)$ equals the distribution of $G(\beta_j)$, which is the $\text{SLE}_\kappa$ trace in $\mathbb{H}$ from $x_j$ to $x_{3-j}$, so we are done. □

**Acknowledgments.** I would like to thank Nikolai Makarov for introducing me to the area of SLE. I also thank Oded Schramm for some important suggestions about this paper and future work.

## REFERENCES

[1] AHLFORS, L. V. (1973). *Conformal Invariants: Topics in Geometric Function Theory*. McGraw-Hill, New York. MR0357743
[2] CAMIA, F. and NEWMAN, C. M. (2007). Critical percolation exploration path and $\text{SLE}_6$: A proof of convergence. *Probab. Theory Related Fields* **139** 473–519. MR2322705
[3] DUBÉDAT, J. (2007). Commutation relations for Schramm–Loewner evolutions. *Comm. Pure Appl. Math.* **60** 1792–1847. MR2358649
[4] LAWLER, G. F., SCHRAMM, O. and WERNER, W. (2001). Values of Brownian intersection exponents. I. Half-plane exponents. *Acta Math.* **187** 237–273. MR1879850
[5] LAWLER, G. F., SCHRAMM, O. and WERNER, W. (2003). Conformal restriction: The chordal case. *J. Amer. Math. Soc.* **16** 917–955. MR1992830
[6] LAWLER, G. F., SCHRAMM, O. and WERNER, W. (2004). Conformal invariance of planar loop-erased random walks and uniform spanning trees. *Ann. Probab.* **32** 939–995. MR2044671
[7] LAWLER, G. F. and WERNER, W. (2004). The Brownian loop soup. *Probab. Theory Related Fields* **128** 565–588. MR2045953
[8] REVUZ, D. and YOR, M. (1991). *Continuous Martingales and Brownian Motion*. Springer, Berlin. MR1083357

[9] ROHDE, S. and SCHRAMM, O. (2005). Basic properties of SLE. *Ann. of Math.* (*2*) **161** 883–924. MR2153402

[10] SCHRAMM, O. (2007). Conformally invariant scaling limits: An overview and a collection of problems. In *International Congress of Mathematicians* **1** 513–543. EMS, Zürich. MR2334202

[11] SCHRAMM, O. (2000). Scaling limits of loop-erased random walks and uniform spanning trees. *Israel J. Math.* **118** 221–288. MR1776084

[12] SCHRAMM, O. and SHEFFIELD, S. Contour lines of the two-dimensional discrete Gaussian free field. Available at arXiv:math.PR/0605337.

[13] SMIRNOV, S. (2001). Critical percolation in the plane: Conformal invariance, Cardy's formula, scaling limits. *C. R. Acad. Sci. Paris Sér. I Math.* **333** 239–244. MR1851632

[14] ZHAN, D. (2004). Stochastic Loewner evolution in doubly connected domains. *Probab. Theory Related Fields* **129** 340–380. MR2128237

[15] ZHAN, D. (2004). Duality of chordal SLE. Available at arXiv:0712.0332. *Invent. Math.* To appear.

[16] ZHAN, D. (2006). Some properties of annulus SLE. *Electron. J. Probab.* **11** 1069–1093. MR2268538

[17] ZHAN, D. (2008). The scaling limits of planar LERW in finitely connected domains. *Ann. Probab.* **36** 467–529.

DEPARTMENT OF MATHEMATICS
YALE UNIVERSITY
P.O. BOX 208283
NEW HAVEN, CONNECTICUT 06520-8283
USA
E-MAIL: dapeng.zhan@yale.edu